\documentclass[12pt]{amsart}
\usepackage{graphicx} 
\usepackage{geometry}
\usepackage[T1]{fontenc}
\usepackage[
  style=alphabetic,
  maxalphanames=99,
]{biblatex}
\renewbibmacro*{volume+number+eid}{%
  \printfield{volume}%
  \printfield{number}%
  \setunit{\addcomma\space}%
  \printfield{eid}}
\DeclareFieldFormat[article]{number}{\mkbibparens{#1}}

\title{Local finiteness in varieties of MS4-algebras}

\author[G. Bezhanishvili]{Guram Bezhanishvili}
\address{Department of Mathematical Sciences\\
New Mexico State University\\
Las Cruces NM 88003\\
USA}
\email{guram@nmsu.edu}

\author[Ch. Meadors]{Chase Meadors}
\address{Department of Mathematics\\
University of Colorado Boulder\\
Boulder, CO 80309\\
USA}
\email{Chase.Meadors@colorado.edu}


\usepackage{mathtools}
\usepackage{amsmath}
\usepackage{amssymb}
\usepackage{mathrsfs}
\usepackage{array}

\usepackage{hyperref}
\hypersetup{
    urlcolor=blue,
}

\usepackage{enumerate}
\usepackage[shortlabels]{enumitem}
\usepackage{multicol}

\usepackage{setspace}

\usepackage{graphicx}

\usepackage{changepage}
\usepackage{wrapfig}

\usepackage{pifont}

\usepackage{cellspace}
\usepackage{listings}

\usepackage[normalem]{ulem}

\usepackage{outlines}

\usepackage{tikz}
\usepackage{tikz-cd}

\usepackage{amsthm}
\counterwithin{figure}{section}

\usepackage[capitalise]{cleveref}

\theoremstyle{definition}
\newtheorem{theorem}{Theorem}[section]
\newtheorem{definition}[theorem]{Definition}
\newtheorem{lemma}[theorem]{Lemma}
\newtheorem{corollary}[theorem]{Corollary}

\newtheorem{question}[theorem]{Question}

\newtheorem{construction}[theorem]{Construction}

\theoremstyle{remark}
\newtheorem{remark}[theorem]{Remark}
\newtheorem{example}[theorem]{Example}

\DeclareMathOperator{\Hom}{Hom}

\DeclareMathOperator{\Clp}{Clp}
\DeclareMathOperator{\Uf}{Uf}

\newcommand{\bb}[1]{\mathbf{#1}}
\newcommand{\frk}[1]{\mathfrak{#1}}
\newcommand{\mca}[1]{\mathcal{#1}}
\newcommand{\msf}[1]{\mathsf{#1}}

\newcommand{\setbuilder}[2]{\left\{#1 : #2\right\}}

\newcommand{\set}[1]{\left\{#1\right\}}

\newcommand{\tand}{\text{ and }}

\newcommand{\la}{\langle}
\newcommand{\ra}{\rangle}

\DeclarePairedDelimiter{\abs}{\lvert}{\rvert}

\newcommand{\what}[1]{\widehat{#1}}

\usepackage{caption}
\usepackage{todonotes}

\newcommand{\msfv}{\mathbf{MS4_S}}
\newcommand{\msfl}{\mathsf{MS4_S}}

\begin{document}

\begin{abstract}
It is a classic result of Segerberg and Maksimova that a variety of $\mathsf{S4}$-algebras is locally finite iff it is of finite depth. 
Since the logic $\mathsf{MS4}$ (monadic $\msf{S4}$) axiomatizes the one-variable fragment of $\mathsf{QS4}$ (predicate $\mathsf{S4}$), it is natural to try to generalize the Segerberg--Maksimova theorem to this setting. 
We obtain several results in this direction. Our positive results include the identification of the largest semisimple variety of $\mathsf{MS4}$-algebras. We prove that the corresponding logic $\msfl$ has the finite model property. We show that both $\msf{S5}^2$ and $\msf{S4}_u$ are proper extensions of $\msfl$, and that a direct generalization of the Segerberg--Maksimova theorem holds for a family of varieties containing the variety of $\msf{S4}_u$-algebras. Our negative results   
include a translation of varieties of $\msf{S5}_2$-algebras into varieties of $\msfl$-algebras of depth 2, which preserves and reflects local finiteness. This, in particular, shows that the problem of characterizing locally finite varieties of $\msf{MS4}$-algebras (even of $\msfl$-algebras) is at least as hard as that of characterizing locally finite varieties of $\msf{S5}_2$-algebras---a problem that remains wide open.
\end{abstract}

\maketitle

\tableofcontents

\section{Introduction}

It is a classic result in modal logic, known as the Segerberg--Maksimova theorem, that a normal extension of $\msf{S4}$ is locally tabular iff it is of finite depth (\cite{Seg71}, \cite{Mak75}). 
Since a logic is locally tabular iff its corresponding variety is locally finite, this provides a characterization of locally finite varieties of $\msf{S4}$-algebras (see \cref{local-finiteness}).

An important extension of $\msf{S4}$ is the bimodal logic $\msf{MS4}$---monadic $\msf{S4}$---which axiomatizes the one-variable fragment of predicate $\msf{S4}$ \cite{FS77}.
It is known that a direct adaptation of the Segerberg--Maksimova theorem is no longer true for varieties of $\msf{MS4}$-algebras. Indeed, the well-known variety of two-dimensional diagonal-free cylindric algebras is not locally finite \cite[Thm.~2.1.11]{HMT85},   
although it is precisely the variety of $\msf{MS4}$-algebras of depth 1. This variety corresponds to the well-known bimodal system $\msf{S5}^2$, which axiomatizes the (diagonal-free) two-variable fragment of classical predicate logic (or, equivalently, the one-variable fragment of predicate $\msf{S5}$). 
Things improve when the attention is restricted to the bimodal logic $\msf{MGrz}$---the one-variable fragment of predicate $\msf{Grz}$ (Grzegorczyk's well-known modal logic \cite{Grz67}). In this case, the Segerberg--Maksimova theorem has a direct generalization: a variety $\bb{V}$ of $\msf{MGrz}$-algebras is locally finite iff it is of finite depth \cite[Sec.~4.10]{Bez01}.

It is natural to seek a characterization of locally finite varieties of $\msf{MS4}$-algebras. In the special case of $\msf{S5}^2$-algebras, such a characterization was given in \cite[Sec.~4]{Bez02}, where it was shown that every proper subvariety of the variety $\bb{S5}^2$ of all $\msf{S5}^2$-algebras is locally finite. 
As we pointed out above, $\msf{S5}^2$-algebras are exactly $\msf{MS4}$-algebras of depth 1. One of our main results (\cref{translation}) shows that already a characterization of locally finite varieties of $\msf{S5}^2$-algebras of depth 2 
is at least as hard as characterizing locally finite varieties of $\bb{S5}_2$---the variety corresponding to the fusion of $\msf{S5}$ with itself.
This problem, as well as the related problem of characterizing locally finite varieties of (unimodal) $\msf{KTB}$-algebras, remains wide open (see \cite{Shap22} for a recent treatment).
This we do by demonstrating a translation of subvarieties of $\bb{S5}_2$ into subvarieties of $\bb{MS4}[2]$---the variety of $\msf{MS4}$-algebras of depth 2---that preserves and reflects local finiteness. Thus, already characterizing locally finite subvarieties of $\bb{MS4}[2]$ would solve the problem of local finiteness for varieties of $\msf{S5}_2$-algebras. 
In addition, we show that our translation does not naturally extend to depth 3 or higher, suggesting that in the general case local finiteness in $\bb{MS4}$ might be even more difficult.
To summarize, 
characterizing locally finite varieties of $\msf{MS4}$-algebras is a hard problem.

On the positive side, we show in \cref{local-finiteness} that the Segerberg--Maksimova theorem has an obvious generalization to a family of varieties containing $\bb{S4}_u$ which corresponds to the well-known bimodal logic $\msf{S4}_u$---$\msf{S4}$ with the universal modality---a logic that plays an important role in the study of the region-based theory of space (see, e.g., \cite{Aie07}). Both $\bb{S5}^2$ and $\bb{S4}_u$ are semisimple subvarieties of $\bb{MS4}$. We identify in \cref{semisimple} the largest semisimple subvariety of $\bb{MS4}$, which we denote by $\msfv$. We prove that the corresponding bimodal logic $\msfl$ has the finite model property. Unfortunately, it is a hard problem to characterize local finiteness already in subvarieties of $\msfv$ as our translation lands in $\msfv[2]$.

These difficulties manifest themselves in the structure of dual spaces of finitely generated $\msf{MS4}$-algebras. While the well-known coloring technique \cite{Esa77} does generalize to $\msf{MS4}$-algebras, unlike finitely generated $\msf{S4}$-algebras, the dual spaces of finitely generated $\msf{MS4}$-algebras may have infinitely many points of finite depth. We demonstrate this by describing the dual space of the well-known Erd\"os--Tarski algebra, which is an infinite one-generated $\msf{S5}^2$-algebra, thus has infinitely many points of depth 1. Some related open problems are posed in \cref{finitely-generated-algebras}, and possible future directions of research are discussed in \cref{sec: conclusion}.

\section{Preliminaries} 

Let $\mathcal L$ be a propositional modal language with two modalities $\lozenge$ and $\exists$. As usual, we write $\square = \neg \lozenge \neg$ and $\forall = \neg \exists \neg$.

\begin{definition}[\cite{FS77}]
    The bimodal logic $\mathsf{MS4}$ is the smallest normal modal logic in $\mathcal L$ containing 
\begin{itemize}
    \item The $\mathsf{S4}$ axioms for $\lozenge$ (i.e.~the $\msf{K}$ axiom along with $p \to \lozenge p$ and $\lozenge \lozenge p \to \lozenge p$),
    \item The $\mathsf{S5}$ axioms for $\exists$ (e.g.~the $\msf{S4}$ axioms along with $\exists \forall p \to \forall p$),
    \item The left commutativity axiom $\exists \lozenge p \to \lozenge \exists p$.
\end{itemize}
\end{definition}

\begin{remark}
    In the terminology of \cite{GKWZ03}, $\msf{MS4} = [\msf{S4}, \msf{S5}]^\text{EX}$; that is, $\msf{MS4}$ is the \textit{expanding relativized product} of $\msf{S4}$ and $\msf{S5}$.
\end{remark}

Algebraic semantics for $\msf{MS4}$ is given by the following Boolean algebras with operators, first considered by Fischer-Servi \cite{FS77} under the name of \textit{bimodal algebras}.

\begin{definition}
\label{def:ms4}
An $\msf{MS4}$-algebra is a tuple $\frk{A} = (B, \lozenge, \exists)$ such that
\begin{itemize}
    \item $(B, \lozenge)$ is an $\msf{S4}$-algebra, i.e.~$B$ is a Boolean algebra and $\lozenge$ is a unary function on $B$ satisfying the identities of a closure operator:
    \[
    \lozenge 0 = 0 \qquad \lozenge (a \vee b) = \lozenge a \vee \lozenge b \qquad a \leq \lozenge a \qquad \lozenge \lozenge a \leq \lozenge a.
    \]
    \item $(B, \exists)$ is an $\msf{S5}$-algebra, that is an $\msf{S4}$-algebra that in addition satisfies $\exists \forall a \leq \forall a$ (where we write $-$ for Boolean negation and use standard abbreviations $\square = {-} \lozenge {-}$ and $\forall = {-} \exists {-}$).
    \item $\exists \lozenge a \leq \lozenge \exists a$.
\end{itemize}
\end{definition}

It is clear that the class of $\msf{MS4}$-algebras is equationally definable, and hence forms a variety. We denote it and the corresponding category by $\bb{MS4}$. 
The last axiom in the preceding definition has a few equivalent forms that we will make use of in the sequel:

\begin{lemma}
    \label{lem:equivalent-axioms}
    Each of the following identities is equivalent to the axiom $\exists \lozenge \leq \lozenge \exists$:
    \begin{enumerate}
        [label=\normalfont(\arabic*), ref = \thelemma(\arabic*)]
        \item \label[lemma]{lem:equivalent-axioms:1} $\exists \lozenge \exists = \lozenge \exists$ ($\exists$ preserves $\lozenge$-fixpoints).
        \item \label[lemma]{lem:equivalent-axioms:2} $\forall \square \forall = \square \forall$ ($\forall$ preserves $\square$-fixpoints).
        \item \label[lemma]{lem:equivalent-axioms:3} $\exists \square \leq \square \exists$.
        \item \label[lemma]{lem:equivalent-axioms:4} $\lozenge \forall \lozenge = \forall \lozenge$ ($\lozenge$ preserves $\forall$-fixpoints).
    \end{enumerate}
\end{lemma}

\begin{proof}
    That the original axiom follows from (1) is obvious since $\exists$ is increasing. 
    For the same reason, the axiom implies the $(\geq)$ direction in (1), while the $(\leq)$ direction follows 
    since $\exists$ is idempotent.
    By taking Boolean negation, we see that (2) is equivalent to (1).
    That the axiom is equivalent to (3) is known, but not immediate;
    see \cite[Lem.~2.5]{Bez23}.
    Taking Boolean negation yields that (3)
    is equivalent to the inequality $\lozenge \forall \leq \forall \lozenge$, which in turn is equivalent to (4)
    in the same manner that (1) is equivalent to the original axiom. 
\end{proof}

\begin{definition}
    \label{def:b0}
    For an $\msf{MS4}$-algebra $\frk{A} = (B, \lozenge, \exists)$ let $B_0$
    be the set of $\exists$-fixpoints.
\end{definition}

If $\lozenge_0$ is the restriction of $\lozenge$ to $B_0$, then it is straightforward to see that $(B_0, \lozenge_0)$ is an $\msf{S4}$-subalgebra of $(B, \lozenge)$. We denote it by $\frk{A}_0$.

Write $H_\square$ 
for the set of $\square$-fixpoints of an $\msf{S4}$-algebra $(B, \lozenge)$. It is well known (see, e.g., \cite[Prop.~2.2.4]{Esa19}) that $H_\square$ is a bounded sublattice of $B$ which is a Heyting algebra, with Heyting implication given by $a \to b = \square (-a \vee b)$.
The Heyting algebra $H_\square$ is sometimes called the \textit{skeleton} of $(B, \lozenge)$, and it is well known that, up to isomorphism, each Heyting algebra arises as the skeleton of some $\msf{S4}$-algebra (see, e.g., \cite[Sec.~2.5]{Esa19}). This plays a key role in proving faithfulness of the well-known G\"{o}del translation of intuitionistic logic into $\msf{S4}$ (see, e.g., \cite[Sec.~XI.8]{RS70}).

\begin{remark}
\label{rem:ms4-s4-representation}
Similar to monadic Heyting algebras (see, e.g., \cite[Sec.~3]{Bez98}), 
we have that $\msf{MS4}$-algebras can be represented as pairs of $\msf{S4}$-algebras $(\frk{A}, \frk{A}_0)$ such that $\frk{A}_0$ is an $\msf{S4}$-subalgebra of $\frk{A}$ and the inclusion $\frk{A}_0 \hookrightarrow \frk{A}$ has a left adjoint $(\exists)$.
\end{remark}

We will make use of J\'{o}nsson--Tarski duality \cite{JT51} to work with the dual spaces of $\msf{MS4}$-algebras. We recall some of the relevant notions:

\begin{definition} 
    A topological space $X$ is a \textit{Stone space} if it is compact, Hausdorff, and zero-dimensional---that is, $X$ has a basis of {\em clopen sets} (sets that are both closed and open).
\end{definition}

\begin{definition} \label{def: cont rel}
    Let $X$ be a Stone space.
    We say that a relation $R \subseteq X^2$ is \textit{continuous} if 
    \begin{enumerate}
        \item $R(x) := \setbuilder{y \in X}{x R y}$ is closed for each $x \in X$ ($R$ is \textit{point-closed}).
        \item $R^{-1}(U) := \setbuilder{y \in X}{yRx \text{ for some } x \in U}$ is clopen whenever $U \subseteq X$ is clopen.
    \end{enumerate}
\end{definition}

Specializing multimodal descriptive general frames (see, e.g., \cite[Sec.~5.5]{BRV01}) to $\msf{MS4}$ yields: 

\begin{definition}
A {\em descriptive $\msf{MS4}$-frame} is a tuple $\frk{F} = (X, R, E)$ such that $X$ is a Stone space, $R$ is a quasi-order (reflexive and transitive), and $E$ is an equivalence relation on $X$ such that
\begin{itemize}
    \label{def:ms4-frame}
    \item both $R$ and $E$ are continuous relations, 
    \item $RE \subseteq ER$---that is, $\forall x,y,y' \in X \; (x E y \tand y R y') \to \exists x' \in X \; (x R x' \tand x' E y')$.

    \begin{center}
    \begin{tikzpicture}[
        scale=1.5,
    	dot/.style={circle,fill=black,minimum size=6pt,inner sep=0}
    ]
    \node[dot,label=left:$x$] (bl) at (0, 0) {};
    \node [dot,label=left:$x'$] (tl) at (0, 1) {};
    \node [dot,label=right:$y'$] (tr) at (1, 1) {};
    \node [dot,label=right:$y$] (br) at (1, 0) {};
    \draw[latex-latex] (bl) -- (br);
    \draw[latex-latex,dashed] (tl) -- (tr);
    \draw[-latex,dashed] (bl) -- (tl);
    \draw[-latex] (br) -- (tr);
    \node at (0.5, -0.2) {$E$};
    \node at (-0.2, 0.5) {$R$};
    \node at (0.5, 1.2) {$E$};
    \node at (1.2, 0.5) {$R$};
    \end{tikzpicture}
    \end{center}
\end{itemize}
\end{definition}

We will refer to descriptive frames simply as \textit{frames} because they are the only kind of frames we will deal with in this paper (as opposed to Kripke frames or the broader class of general frames studied in modal logic \cite{CZ97}, \cite{BRV01}).

\begin{remark} \label{rem: subalgebras}
    The subalgebra $\frk{A}_0 = (B_0, \lozenge_0)$ of $\exists$-fixpoints (\cref{def:b0}) is dually identified with the $\msf{S4}$-frame $(X/E, \overline{R})$ of $E$-equivalence classes, where 
    \[
    \alpha \overline{R} \beta \ \mbox{ iff } \ \exists x \in \alpha, y \in \beta : xRy.
    \] 
    The condition $RE \subseteq ER$ ensures that this is a well-defined quasi-order.
    In fact, this is a general construction arising from any \textit{correct partition} of a descriptive frame (\cite{Esa77}). We discuss this in more detail in \cref{finitely-generated-algebras}.
\end{remark}

The next definition is well known (see, e.g., \cite[Sec.~3.3]{BRV01}).

\begin{definition} 
    Let $X, X'$ be sets, $S$ a binary relation on $X$, and $S'$ a binary relation on $X'$.
    A function $f : X \to X'$ is a \textit{p-morphism} (or \textit{bounded morphism})
    with respect to $(S, S')$ if 
    \begin{itemize}
        \item $x S y$ implies $f(x) S' f(y)$;
        \item $f(x) S' y'$ implies $x S y$ for some $y \in X$ with $f(y) = y'$.
    \end{itemize}
\end{definition}

\begin{definition}
    \label{def:ms4-p-morphism}
    An {\em $\msf{MS4}$-morphism} is a continuous map $f:(X, R, E) \to (X', R', E')$ between $\msf{MS4}$-frames such that $f$ is a p-morphism with respect to both $(R, R')$ and $(E, E')$.
\end{definition}

Specializing J\'{o}nsson--Tarski duality 
to $\msf{MS4}$-algebras yields:

\begin{theorem}
    The category of $\msf{MS4}$-algebras with homomorphisms and the category of descriptive $\msf{MS4}$-frames with $\msf{MS4}$-morphisms are dually equivalent.
\end{theorem}

\begin{remark}
\label{rem:jt-duality}
The above dual equivalence is implemented 
as follows:
with each $\msf{MS4}$-frame $\mathfrak{F} = (X, R, E)$ we associate the $\msf{MS4}$-algebra $\mathfrak{F}^* = (\Clp X, \lozenge_R, \lozenge_E)$ of clopen subsets of $X$, where $\lozenge_R$ and $\lozenge_E$ are the \textit{dual operators} of $R$ and $E$, defined by
\[
\lozenge_R U = R^{-1}(U) \ \mbox{ and } \ \lozenge_E U = E(U).
\]
In the other direction, with each $\msf{MS4}$-algebra $\frk{A} = (B, \lozenge, \exists)$ we associate the descriptive frame $\frk{A}_* = (\Uf B, R_\lozenge, R_\exists)$ of ultrafilters of $B$, where $R_\lozenge$ and $R_\exists$ are the \textit{dual relations} of $\lozenge$ and $\exists$, defined by 
\[
x R_\lozenge y \mbox{ iff } x \cap H_\square \subseteq y \ \mbox{ and } \ x R_\exists y \mbox{ iff } x \cap \frk{A}_0 = y \cap \frk{A}_0,
\]
and $\Uf B$ is equipped with the \emph{Stone topology} generated by the clopen basis $\setbuilder{\frk{s}(a)}{a \in B}$, where $\frk{s}(a) = \setbuilder{x}{a \in x}$.
The unit isomorphisms of the duality are given by
\begin{itemize}
    \item $\frk{A} \to (\frk{A}_*)^*$; \ $a \mapsto \setbuilder{x}{a \in x}$,
    \item $\frk{F} \to (\frk{F}^*)_*$; \ $x \mapsto \setbuilder{a}{x \in a}$.
\end{itemize}
The dual of a morphism $f$ (in both categories) is the inverse image map $f^{-1}[\cdot]$, and 
the natural bijection $\Hom(\frk{A}, \frk{F}^*) \cong \Hom(\frk{F}, \frk{A}_*)$ is given by associating to $f : \frk{A} \to \frk{F}^*$ the morphism $\bar{f} : \frk{F} \to \frk{A}_*$ defined by $x \mapsto \setbuilder{a \in \frk{A}}{x \in f(a)}$, and to $g : \frk{F} \to \frk{A}_*$ the morphism $\bar{g} : \frk{A} \to \frk{F}^*$ defined by $a \mapsto \setbuilder{x \in \frk{F}}{a \in g(x)}$.
\end{remark}

\section{Semisimple varieties}
\label{semisimple}

Recall (\cite[Sec.~IV.12]{BSS81}) that a variety $\bb{V}$ of algebras is \textit{semisimple} if every subdirectly irreducible algebra from $\bb{V}$ is simple.
It is well known
that $\bb{S5}$ is the largest semisimple subvariety of $\bb{S4}$. 
In this section we introduce the largest semisimple subvariety of $\bb{MS4}$ and its corresponding logic, which we denote by $\msfv$ and $\msfl$, respectively.
The variety $\msfv$ contains two well-known varieties 
$\bb{S5}^2$ and $\bb{S4}_u$, whose corresponding logics are $\msf{S5}^2$ (\cite[230]{GKWZ03}) and $\msf{S4}_u$ (\cite[38]{GKWZ03}).
One of our main results in this section shows that $\msfv$ is generated by its finite algebras, and hence that $\msf{MS4_S}$ has the finite model property (and is decidable).

The standard correspondence between congruences and modal filters in BAOs (see \cite[Sec.~4]{Ven07}) yields that the lattice of congruences of an $\bb{MS4}$-algebra $\frk{A}$ is isomorphic to the lattice of modal filters of $\frk{A}$, where  a \textit{modal filter} of $\frk{A}$ is a filter $F$ satisfying $a \in F$ implies $\square a, \forall a \in F$.
Under standard duality theory (as outlined in \cref{rem:jt-duality}), filters of $\frk{A}$ closed under a modal necessity operator $\square$ correspond dually to closed subsets of $\frk{A}_*$ that are $R_\lozenge$-\textit{upsets}, i.e.~closed subsets $U$ such that $x \in U$ and $x R_\lozenge y$ imply $y \in U$ (where $R_\lozenge$ is the dual relation of $\lozenge$).
Hence, modal filters correspond dually to closed $E$-saturated $R$-upsets (we refer to $E$-upsets as \textit{$E$-saturated sets} since they are unions of equivalence classes of $E$).
We next give a nicer characterization of congruences of $\frk{A}$. For this we introduce an auxiliary $\msf{S4}$-modality on $\frk{A}$, and the corresponding auxiliary quasi-order on the dual descriptive frame of $\frk{A}$.

\begin{definition}\
\begin{enumerate}
    \item For an $\msf{MS4}$-algebra $\mathfrak A=(B, \lozenge, \exists)$, define $\blacklozenge = \lozenge \exists$.
    \item For an $\msf{MS4}$-frame $\mathfrak F=(X, R, E)$, define $Q = ER$.
\end{enumerate}
\end{definition}

Clearly $\blacklozenge$ is an $\bb{S4}$-possibility operator on $B$ and $Q$ is a quasi-order on $X$. The next lemma shows that $Q$ is the dual relation of $\blacklozenge$.
Let $\blacksquare = \neg \blacklozenge \neg$. Then $\blacksquare = \Box\forall$ and $\blacksquare$ is an $\msf{S4}$-necessity operator on $B$. We let $H_\blacksquare$ be the Heyting algebra of fixpoints of $\blacksquare$.

\begin{lemma}
    Let $\mathfrak A$ be an $\msf{MS4}$-algebra and $\mathfrak F$ its dual $\msf{MS4}$-frame. Then $xQy$ iff $x\cap H_\blacksquare\subseteq y$ for each $x,y\in \mathfrak F$. 
\end{lemma}

\begin{proof}
    First suppose that $xQy$ and $a\in x\cap H_\blacksquare$. Then there is $z$ such that $xRz$ and $zEy$. Therefore, $x\cap H_\square \subseteq z$ and $z\cap B_0 =y\cap B_0$. Since $a\in H_\blacksquare$, we have $a=\Box\forall a=\forall\Box\forall a$ (where the last equality follows from \cref{lem:equivalent-axioms:2}), so $a\in H_\square\cap B_0$. 
    But then from $a\in x$ it follows that $a\in z$, which then implies that $a\in y$. Thus, $x\cap H_\blacksquare \subseteq y$.
    
    Conversely, suppose that $x\cap H_\blacksquare \subseteq y$. Let $F$ be the filter generated by $(x\cap H_\square)\cup(y\cap B_0)$ and $I$ the ideal generated by $B_0\setminus y$. We show that $F \cap I = \varnothing$. Otherwise there are $a\in x\cap H_\square$, $b\in y\cap B_0$, and $c\in B_0\setminus y$ such that $a\wedge b \le c$. Therefore, $a\le b\to c$, so $a \le \Box (b \to c)$ because $a\in H_\square$. Since $b,c\in B_0$, we have $\Box (b\to c)\in H_\blacksquare$. Thus, $\Box(b\to c)\in x\cap H_\blacksquare\subseteq y$, and so $b\to c\in y$. Consequently, $c \in y$ because $b \in y$. The obtained contradiction proves that $F\cap I=\varnothing$.
    Therefore, there is an ultrafilter $z$ such that $x\cap H_\square\subseteq z$ and $z\cap B_0=y\cap B_0$. But then $xRz$ and $zEy$, yielding that $xQy$.
\end{proof}

\begin{definition}
 We call a filter $F$ of an $\msf{MS4}$-algebra a \textit{$\blacksquare$-filter} if $a\in F$ implies $\blacksquare a \in F$.
\end{definition}

We order the set of $\blacksquare$-filters by inclusion. It is then clear that it is a complete lattice (because it is closed under arbitrary intersections).
The next theorem generalizes the well-known correspondence (see, e.g., \cite[Secs.~2.4, 3.4]{Esa19}) between congruences and $\square$-filters of $\msf{S4}$-algebras and closed upsets of their dual $\msf{S4}$-frames to the setting of $\msf{MS4}$-algebras.

\begin{theorem}
    \label{thm:ms4-congruence}
    Let $\frk{A}$ be an $\msf{MS4}$-algebra and $\frk{F}$ its dual $\sf MS4$-frame. The following complete lattices are isomorphic:
    \begin{enumerate}
        \item Modal filters (and congruences) of $\frk{A}$.
        \item $\blacksquare$-filters of $\frk{A}$.
        \item Modal filters (and congruences) of the $\msf{S4}$-algebra $\frk{A}_0$.
        \item Filters (and congruences) of the Heyting algebra $H_\blacksquare$.
        \item Closed $Q$-upsets of $\frk{F}$ (under reverse inclusion).
        \item Closed $E$-saturated $R$-upsets of $\frk{F}$ (under reverse inclusion).
    \end{enumerate}
\end{theorem}

\begin{proof}
The modal filters and $\blacksquare$-filters are literally the same lattice. Therefore, so are the closed $Q$-upsets and the closed $E$-saturated $R$-upsets of $\frk{F}$.
This gives the equivalence of (1), (2) and (5), (6). 
Note that the $\blacksquare$-fixpoints of $\frk{A}$ are exactly the $\square$-fixpoints of $\frk{A}_0$, which yields the equivalence of (2) and (3), and also (4) by 
\cite[Thm.~2.4.17]{Esa19}. Finally, these are equivalent to (5) by \cite[Thm.~3.4.16]{Esa19}.
\end{proof}

Let $\frk{F} = (X, R, E)$ be an $\sf MS4$-frame. We call $x \in X$ a \textit{$Q$-root} of $\frk{F}$ if $Q(x) = X$.
Then $\frk{F}$ is \textit{$Q$-rooted} if it has a $Q$-root, and \textit{strongly $Q$-rooted} if the set of $Q$-roots is nonempty and open.
A \textit{$Q$-cluster} is an equivalence class of $Q \cap Q^{-1}$.
The previous theorem then yields a characterization of simple and subdirectly irreducible algebras (we abbreviate subdirectly irreducible by s.i.):

\begin{theorem}
    \label{thm:si-simple-ms4}
    Let $\frk{A}$ be a (non-trivial) $\msf{MS4}$-algebra and $\frk{F}$ its dual $\sf MS4$-frame.
    \begin{enumerate}
        [label=\normalfont(\arabic*), ref = \thetheorem(\arabic*)]
        \item \label[theorem]{thm:si-simple-ms4:1} $\frk{A}$ is s.i.~iff $\frk{F}$ is strongly $Q$-rooted.
        \item \label[theorem]{thm:si-simple-ms4:2} $\frk{A}$ is simple iff $\frk{F}$ is a $Q$-cluster.
    \end{enumerate}
\end{theorem}

\begin{proof}
    (1) It is well known (see, e.g. \cite[Sec.~II.8]{BSS81}) that $\frk{A}$ is s.i.~iff $\frk{A}$ has a least nontrivial congruence.
    By the above theorem, this is equivalent to $\frk{F}$ having the largest closed $Q$-upset $Y$ different from $X$. The complement of $Y$ is the nonempty open set of $Q$-roots of $\frk{F}$. 

    (2) $\frk{A}$ is simple iff the only nontrivial congruence is $A^2$ (\cite[Sec.~II.8]{BSS81}).
    Thus, using (1), $\frk{A}$ is simple iff $Y=\varnothing$, which is equivalent to $\frk{F}$ being a $Q$-cluster. 
\end{proof}

\begin{remark}
    The above theorem also follows from a general result of 
    Venema \cite{Ven04}. For this it is enough to observe that $Q$ is the reachability relation in an $\msf{MS4}$-frame, and thus the above theorem follows from \cite[Cor.~1]{Ven04}.
\end{remark}

Let $\frk{A}$ be a s.i.~$\msf{MS4}$-algebra. By the above theorem, $\frk{A}$ is simple iff $\blacklozenge$ is an $\msf{S5}$-possibility operator, which happens iff $\blacklozenge\blacksquare a \le \blacksquare a$ for each $a\in \frk{A}$. This motivates the following definition.

\begin{definition}
    Let $\msf{MS4_S}=\msf{MS4}+\blacklozenge \blacksquare p \to \blacksquare p$ and let $\msfv = \bb{MS4}+(\blacklozenge \blacksquare a \leq \blacksquare a)$ be the corresponding variety. 
\end{definition}

\begin{theorem}
    \label{thm:ms4s-characterization}
    $\msfv$ is the largest semisimple subvariety of $\bb{MS4}$.
\end{theorem}

\begin{proof}
Let $\bb{V}$ be a semisimple subvariety of $\bb{MS4}$ and let $\frk{A} \in \bb{V}$ be s.i.
Then $\frk{A}$ is simple, so the dual $\sf MS4$-frame $\frk{F}$ of $\frk{A}$ is a $Q$-cluster by \cref{thm:si-simple-ms4:2}.
Therefore, $Q$ is an equivalence relation on $\frk{F}$, and hence $\frk{A}\models\blacklozenge \blacksquare p \to \blacksquare p$.
Thus, $\frk{A}\in\msfv$, and so $\bb{V}\subseteq\msfv$.
\end{proof}

\begin{remark}
    The above theorem also follows from a general result of Kowalski and Kracht   \cite[Thm.~12, Prop.~8]{Kow06} that semisimplicity of a variety of modal algebras in finite signature is equivalent 
    to the definability of a universal modality.
\end{remark}

The two logics $\msf{S4}_u$ and $\msf{S5}^2$ mentioned in the introduction are actually both extensions of $\msfl$.
We recall the definitions:

\begin{definition}\
\begin{enumerate} 
    \item \cite{GVP92,Ben96} $\msf{S4}_u$ is the smallest normal modal logic in $\mathcal L$
    containing the $\msf{S4}$ axioms for $\lozenge$, the $\msf{S5}$ axioms for $\exists$, and the bridge axiom $\lozenge p \to \exists p$.
    Algebraically, $\bb{S4}_u$ is the variety $\bb{MS4} + (\lozenge a \leq \exists a)$.
    \item \cite[Sec.~5.1]{GKWZ03}
    $\msf{S5}^2$ is the \textit{product} of $\msf{S5}$ with itself, i.e.~the smallest normal modal logic in $\mathcal L$
    containing the $\msf{S5}$ axioms for both $\lozenge,\exists$ and the commutativity axiom $\lozenge \exists p \leftrightarrow \exists \lozenge p$. 
    Algebraically, $\bb{S5}^2$ is the variety $\bb{MS4} + (\lozenge \square a \leq \square a)$.
\end{enumerate}
\end{definition}

The variety $\bb{S5}^2$ is also known as the variety of \textit{diagonal-free cylindric algebras of dimension two} \cite[Def.~1.1.2]{HMT85}. Its subvarieties were studied in detail in \cite{Bez02}. 
As we will see in \cref{local-finiteness}, $\bb{S5}^2$ is the variety of all $\msf{MS4}$-algebras of \textit{depth 1}.

\begin{figure}[h]
\centering
\begin{tikzpicture}[
    scale=1.25,
	dot/.style={circle,fill=black,minimum size=6pt,inner sep=0}
]
\node[dot] (top) at (0, 1) {};
\draw[-latex] (0, 0) node[dot] {} -- (top);
\draw[dotted] (-0.2, -0.2) rectangle (0.2, 1.2);
\node at (0, -0.5) {$\frk{F}$};
\end{tikzpicture}
\qquad 
\begin{tikzpicture}[
    scale=1.25,
	dot/.style={circle,fill=black,minimum size=6pt,inner sep=0}
]
\draw (0, 0) node[dot] {} -- (1, 0) node[dot] {};
\draw (0, 1) node[dot] {} -- (1, 1) node[dot] {};
\draw[dotted] (-0.2, -0.2) rectangle (0.2, 1.2);
\draw[dotted] (0.8, -0.2) rectangle (1.2, 1.2);
\node at (0.5, -0.5) {$\frk{G}$};
\end{tikzpicture}
\caption{
    The dual algebra of $\frk{F}$ 
    belongs to $\bb{S4}_u$ but not to $\bb{S5}^2$, and the opposite holds for the dual algebra of $\frk{G}$.
}
\label{fig:s52-s4u-sep}
\end{figure}
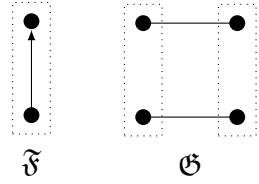

Note that $\bb{S4}_u$ and $\bb{S5}^2$ are incomparable. 
Indeed, consider the frames $\frk{F}$ and $\frk{G}$ depicted in \cref{fig:s52-s4u-sep}, where the black arrows represent $R$ and dotted circles $E$. The dual algebra of $\frk{F}$ belongs to $\bb{S4}_u$ since there is only one $E$-cluster, but evidently $R$ is not an equivalence, so it does not belong to $\bb{S5}^2$.
On the other hand, the dual algebra of $\frk{G}$ belongs to $\bb{S5}^2$ since $R$ and $E$ are commuting equivalence relations, but $R$-clusters are not contained in $E$-clusters, so the identity $\lozenge a \leq \exists a$ of $\bb{S4}_u$ is falsified.
Consequently, both $\bb{S4}_u$ and $\bb{S5}^2$ are proper subvarieties of $\msfv$. 

While the characterization of locally finite subvarieties of these two varieties is known (see \cite{Bez02} and \cref{cor:s4u-lf-finite-depth} below), \color{black} the question is dramatically more complicated for subvarieties of $\msfv$, as will be demonstrated in \cref{translation}.

We conclude this section by establishing that $\msfv$ is generated by its finite members, and hence that $\msf{MS4_S}$ has the finite model property (hereafter abbreviated \textit{fmp}).
For this we recall (see, e.g., \cite[Def.~II.10.14]{BSS81}) that an algebra is \textit{locally finite} if every finitely generated subalgebra is finite.

\begin{lemma} \label{lem: finite S52}
    If $\frk{A} = (B, \lozenge, \exists)$ is an $\msf{MS4_S}$-algebra, then $(B, \blacklozenge, \exists)$ is a locally finite $\msf{S5}^2$-algebra.
\end{lemma}

\begin{proof}
    Since $\frk{A}$ is an $\msf{MS4_S}$-algebra, $\blacklozenge$ is an $\msf{S5}$-operator. Moreover, 
    \[
    \exists \blacklozenge a = \exists \Diamond \exists a = \Diamond \exists a = \blacklozenge a = \blacklozenge \exists a
    \]
    (the second equality follows from \cref{lem:equivalent-axioms:1}).
    Therefore, $(B, \blacklozenge, \exists)$ is an $\msf{S5}^2$-algebra that in addition satisfies $\exists \blacklozenge a = \blacklozenge a$ for each $a\in B$. Thus, $(B, \blacklozenge, \exists)$ belongs to a proper subvariety of $\bb{S5}^2$, which is locally finite by \cite[Sec.~4]{Bez02}. Consequently, $(B, \blacklozenge, \exists)$ is locally finite.
\end{proof}

Let $\frk{A} = (B, \lozenge, \exists)$ be an $\msf{MS4_S}$-algebra. For a finite $S \subseteq B$,
let $B'$ be the subalgebra of $(B, \blacklozenge, \exists)$ generated by $S$. By the above lemma, $B'$ is finite. 
Let $K$ be the fixpoints of $\lozenge$, define $\lozenge'$ on $B'$ by
\[
\lozenge' a = \bigwedge \{ x \in B' \cap K : a \leq x \},
\]
and set $\frk{A}_S=(B', \lozenge', \exists)$. 

\begin{lemma}
    $\frk{A}_S$ is a finite $\msf{MS4_S}$-algebra.
\end{lemma}

\begin{proof}
    As we already pointed out, $\frk{A}_S$ is finite by \cref{lem: finite S52}.
    Moreover, $\lozenge'$ is an $\bb{S4}$-operator on $\frk{A}_S$ (see \cite[Lem.~2.3]{MT44}). For the reader's convenience, we sketch the proof of this result.
    \begin{itemize}
        \item That $\lozenge' 0 = 0$ is clear since $0 \in B' \cap K$.
        \item That $\lozenge'(a \vee b) = \lozenge' a \vee \lozenge' b$ follows from distributivity since the meets involved in the definition of $\lozenge'$ are finite. 
        \item That $a \leq \lozenge'a$ follows from the definition of $\lozenge' a$. 
        \item To see that $\lozenge' \lozenge' a \leq \lozenge' a$, it is enough to observe that if $b\in B'\cap K$, then $a\le b$ iff $\Diamond' a\le b$.
    \end{itemize}
    We next show that $\exists \lozenge' a \leq \lozenge' \exists a$, 
    or equivalently that $\exists \lozenge' \exists a = \lozenge' \exists a$, i.e.~that $\lozenge' \exists a$ is a $\exists$-fixpoint for each $a\in\frk{A}_S$.
    First observe that if $b \in B' \cap K$ with $\exists a \leq b$, then $\forall b \in B' \cap K$ (since $\lozenge \forall \lozenge = \forall \lozenge$) and $\exists a \leq \forall b$.
    Because $\forall b \leq b$, we have $\lozenge' \exists a = \bigwedge \setbuilder{\forall b}{b \in B' \cap K, \exists a \leq b}$.
    But this is a finite meet of $\exists$-fixpoints and hence
    an $\exists$-fixpoint.

    Finally, we check $\blacklozenge' \blacksquare' a \leq \blacksquare' a$ for each $a\in B'$, where $\blacklozenge' = \lozenge' \exists$ and $\blacksquare' = {-} \blacklozenge' {-}$.
    Since $\lozenge \exists a \in B'$ (because $B'$ is closed under $\blacklozenge$), $\lozenge' \exists a = \lozenge \exists a$, so 
    $\blacklozenge' = \blacklozenge \vert_{B'}$, and hence the inequality holds because $\frk{A}$ is an $\msf{MS4_S}$-algebra.
\end{proof}

\begin{theorem}
    $\msf{MS4_S}$ has the fmp.
\end{theorem}

\begin{proof}
    If $\msfl\not\vdash\varphi$, then there is an $\msf{MS4_S}$-algebra $\frk{A} = (B, \lozenge, \exists)$ and an $n$-tuple $\overline{a}$ such that $\varphi(\overline{a}) \neq 1$ in $\frk{A}$.
    Let $S$ be the set of subterms of $\varphi(\overline{a})$.
    Since $S$ is finite, $\frk{A}_S = (B', \lozenge', \exists)$ is a finite $\msf{MS4_S}$-algebra by the previous lemma.
    Moreover, the definition of $\lozenge'$ ensures that if $\lozenge a \in B'$, then $\lozenge' a = \lozenge a$.
    Therefore, the computation of $\varphi(\overline{a})$ in $\frk{A}_S$ is identical to that 
    in $\frk{A}$, and hence $\varphi(\overline{a}) \neq 1$ in $\frk{A}_S$.
    Thus, $\varphi$ is falsified in $\frk{A}_S$. 
\end{proof}

As a corollary we obtain the following pair of results, where the first concerns the finite-depth extensions of $\msf{MS4}$ (a notion we will define in \cref{local-finiteness}). 

\begin{corollary}\
    \begin{enumerate}
        \item $\msfl = \bigcap_n \msfl[n]$. 
        \item $\msfl$ is decidable.
    \end{enumerate}
\end{corollary}

\begin{proof}
    (1) 
    If $\msfl \not\vdash \varphi$, then by the previous theorem, $\varphi$ is falsified on a finite $\msfl$-algebra $\frk{A}$. 
    Let $n$ be the depth of $\frk{A}$. Then 
    $\frk{A}\models\msfl[n]$, so $\msfl[n] \not\vdash \varphi$, and hence
    $\bigcap_n \msfl[n]\not\vdash\varphi$.

    (2) This is obvious since $\msfl$ is finitely axiomatizable and has the fmp.
\end{proof}

\section{Local finiteness}
\label{local-finiteness}

A classic result of Segerberg and Maksimova 
gives a characterization of local finiteness in the lattice of subvarieties of $\bb{S4}$.
By the \textit{depth} of an $\msf{S4}$-algebra $\frk{A}$ or its dual frame $\frk{F}$, we mean the longest length of a proper $R$-chain in $\frk{F}$;
that is, a sequence of points $x_1 R x_2 R \dots R x_n$ where $\neg (x_j R x_i)$ for any $i < j$.
If there is no bound on such $R$-chains, then $\frk{A}$ is of depth $\omega$.
We say a variety $\bb{V} \subseteq \bb{S4}$ has {\em depth} $\leq n$ if the depth of each algebra in $\bb{V}$ is $\leq n$; if $\bb{V}$ contains algebras of arbitrary depth, we say it is of {\em depth} $\omega$.
We will also apply these notions to (varieties of) $\msf{MS4}$-algebras and their dual frames, understanding the depth of an $\msf{MS4}$-algebra to be the depth of its $\msf{S4}$-reduct (that is, the $R$-depth of its dual frame).

Consider the family of formulas $P_n$ defined by
\[
P_1 = \lozenge \square q_1 \to \square q_1 \qquad P_n = \lozenge (\square q_n \wedge \neg P_{n-1}) \to \square q_n.
\]
It is well known (see, e.g., \cite[Thm.~3.44]{CZ97}) that a variety $\bb{V} \subseteq \bb{S4}$ has depth $\leq n$ iff $\bb{V} \models P_n$.
Following \cite{Shap16}, we write $\bb{V}[n]$ to refer to the subvariety $\bb{V} + P_n$.
For example, since equivalence relations are precisely quasi-orders of depth 1, we have $\bb{S4}[1] = \bb{S5}$.
Similarly, $\bb{S5}^2 = \bb{MS4}[1]$.
There is a pleasant analogy between $\bb{S5} \subseteq \bb{S4}$ and $\msfv \subseteq \bb{MS4}$ in that $\bb{S5}$ is the largest subvariety of $\bb{S4}$ of $R$-depth 1, while $\msfv$ is the largest subvariety of $\bb{MS4}$ of $Q$-depth 1.
We have (see, e.g., \cite[Sec.~12.4]{CZ97}):

\begin{theorem}[Segerberg--Maksimova]
A variety $\bb{V} \subseteq \bb{S4}$ is locally finite iff $\bb{V} \models P_n$ for some $n < \omega$.
\end{theorem}

The lattice of subvarieties of $\bb{S4}$ has 
the structure depicted in \cref{s4-subvariety-lattice}:

\begin{figure}[h]
\centering
\begin{tikzcd}
{} \arrow[rrrrdd, no head, dashed, shift right=2] & \mathbf{S4} \arrow[rrd, no head, bend left=15] \arrow[rrd, no head, bend right=15]                          &  &                                        &    \\
& \vdots \arrow[u, no head]                                                                             &  & \mathbf{Grz.3}                         &    \\
& \mathbf{S4}[3] \arrow[rrd, no head, bend left=15] \arrow[rrd, no head, bend right=15] \arrow[u, no head] &  & \vdots \arrow[u, no head]              & {} \\
& \mathbf{S4}[2] \arrow[rrd, no head, bend left=15] \arrow[rrd, no head, bend right=15] \arrow[u, no head] &  & {\mathbf{Grz.3}[3]} \arrow[u, no head] &    \\
& \mathbf{S5} \arrow[rrd, no head] \arrow[u, no head]                                                   &  & {\mathbf{Grz.3}[2]} \arrow[u, no head] &    \\
&                                                                                                       &  & {\mathbf{Grz.3}[1]} \arrow[u, no head] &   
\end{tikzcd}
\caption{The lattice of subvarieties of $\bb{S4}$, dually isomorphic to the lattice of normal extensions of $\msf{S4}$. The locally finite varieties are precisely the ones below the dotted line.}
\label{s4-subvariety-lattice}
\end{figure}
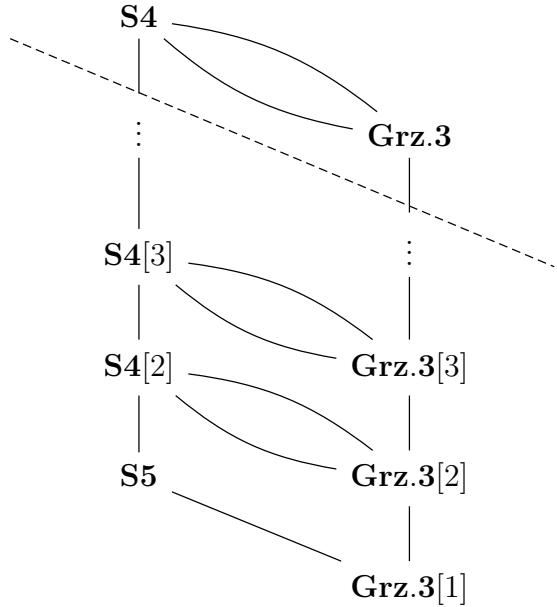

In particular, the variety $\bb{Grz.3}$, the subvariety of $\bb{S4}$ defined by the formulas
\[
\square(\square(p \to \square p) \to p) \to p
\qquad
\square(\square p\to q) \vee \square(\square q\to p),
\]
is the least non-locally-finite subvariety of $\bb{S4}$ (and hence local finiteness below $\bb{S4}$ is decidable).
In the figure, the varieties $\bb{Grz.3}[n]$ appearing on the right are the varieties generated by the single algebra whose dual space is the $n$-element chain, and $\bb{Grz.3}$ is the variety generated by all such algebras.

It is natural to try to extend the Segerberg--Maksimova theorem to $\bb{MS4}$.
However, the theorem fails even in depth $1$ since the variety $\bb{MS4}[1] = \bb{S5}^2$ is not locally finite: the well-known Erd\"os--Tarski algebra 
is a one-generated infinite $\bb{S5}^2$-algebra (see, e.g., \cite[Thm.~2.1.11]{HMT85}). 
Thus, finite depth does not characterize local finiteness in subvarieties of $\bb{MS4}$.
To the authors' knowledge, the only positive result in this direction is that the Segerberg--Maksimova theorem holds for the variety $\bb{MGrz}$ 
of monadic Grzegorczyk algebras (see \cite[Sec.~4.10]{Bez01}).

In this paper we will make the case that characterizing local finiteness in subvarieties of $\bb{MS4}$ is a hard problem.
In fact, we will show that a characterization of local finiteness even in subvarieties of $\msfv[2]$ would yield one for $\bb{S5}_2$.
Here $\bb{S5}_2$ corresponds to the fusion $\msf{S5}_2 = \msf{S5} \oplus \msf{S5}$, the bimodal logic of two $\msf{S5}$ modalities with no connecting axioms (discussed in more detail in \cref{translation}).

We do have the following criterion of local finiteness in subvarieties of $\bb{MS4}$, but item (2) is in general hard to verify (as we will see shortly).

\begin{theorem}
    \label{thm:ms4-lf}
    A subvariety $\bb{V}$ of $\bb{MS4}$ is locally finite iff 
    \begin{enumerate}
        [label=\normalfont(\arabic*), ref = \thetheorem(\arabic*)]
        \item \label[theorem]{thm:ms4-lf:1} $\bb{V}$ is of finite depth ($\bb{V} \models P_k$ for some $k$);
        \item \label[theorem]{thm:ms4-lf:2} there is a function $f:\omega\to\omega$ such that for every $n$-generated s.i.~$\frk{A} \in \bb{V}$, the algebra  $\frk{A}_0$ is $f(n)$-generated as an $\msf{S4}$-algebra.
    \end{enumerate}
\end{theorem}

\begin{proof}
    First suppose that $\bb{V}$ is locally finite. Then for every $n$, the free $n$-generated $\bb{V}$-algebra is finite. Letting $f(n)$ be its cardinality, $f:\omega\to\omega$ is our desired function, and hence (2) is satisfied. Moreover, (1) is satisfied since the $\bb{S4}$-reduct of $\bb{V}$ is a uniformly locally finite class (that is, there is a uniform finite bound on the size of an $n$-generated subalgebra of an algebra from this class), and so it generates a locally finite variety (see \cite[285]{Mal73}). 
    Thus, $\bb{V}$ is of finite depth by the Segerberg--Maksimova theorem. 
    
    Conversely, suppose that $\bb{V}$ satisfies (1) and (2). 
    By \cite[Thm.~3.7(4)]{Bez01}, it suffices to give a finite bound on the cardinality of every $n$-generated s.i.~$\frk{A} \in \bb{V}$.
    By (2), $\frk{A}_0$ is a $f(n)$-generated $\msf{S4}$-algebra, and it is of finite depth because $\frk{A}$ is of finite depth by (1). Thus,  $\abs{B_0} \leq k(n)$ by the Segerberg--Maksimova theorem.
    Since $B_0$ is the $\exists$-fixpoints of $\frk{A}$, we have that $\frk{A}$ is $(n + k(n))$-generated as an $\msf{S4}$-algebra, and thus we have some $m(n)$ such that $\abs{B} \leq m(n + k(n))$.
\end{proof}

We next give a few applications of this theorem where (2) can be verified.
For example, one can verify (2) by induction in the $\bb{MGrz}$ setting, as in \cite[Sec.~4.10]{Bez01}.
For another example, we need the following lemma:

\begin{lemma}
    \label{lem:altk-b0}
    Let $\frk{A}$ be a s.i.~$\msf{MS4}$-algebra. Then, for any $k \geq 1$, $\frk{A}$ validates
    \[
    \msf{alt}_k^0 = \square \forall p_1 \vee \square (\forall p_1 \to \forall p_2) \vee \dots \vee \square(\forall p_1 \wedge \forall p_2 \wedge \dots \wedge \forall p_k \to \forall p_{k+1})    
    \]
    iff $\abs{\frk{A}_0} \leq 2^k$.
\end{lemma}

\begin{proof}
    The formula $\msf{alt}_k^0$ is obtained from the well-known formula 
    \[
    \msf{alt}_k = \square p_1 \vee \square(p_1 \to p_2) \vee \dots \vee \square(p_1 \wedge p_2 \wedge \dots \wedge p_k \to p_{k+1})
    \]
    by substituting $\forall p_i$ for $p_i$. 
    Let $\frk{F} = (X, R, E)$ be the dual of $\frk{A}$. Since $\frk{A}$ is s.i., $\frk{F}$ is strongly $Q$-rooted by \cref{thm:si-simple-ms4:1}. 
    Because the universe of the subalgebra $\frk{A}_0$ is exactly $\forall$-fixpoints of $\frk{A}$, $\frk{A} \models \msf{alt}_k^0$ iff $\frk{A}_0 \models \msf{alt}_k$.
    By \cref{rem: subalgebras}, the dual space of $\frk{A}_0$ is $\frk{F}_0 = (X/E, \overline{R})$, which is a rooted descriptive $\msf{S4}$-frame.
    Since $\msf{alt}_k$ is a canonical formula (\cite[Def.~4.30]{BRV01} together with \cite[Thm.~5.16]{CZ97}), $\msf{alt}_k$ is valid on $\frk{F}_0$ iff it is valid on the underlying Kripke frame of $\frk{F}_0$ (\cite[Prop.~5.85]{BRV01}).
    It is well known (\cite[Prop.~3.45]{CZ97}) that $\msf{alt}_k$ holds in a Kripke frame iff the number of alternatives of each point is $\leq k$ or, in the rooted transitive case, the number of points in the frame is $\leq k$. Thus, $\msf{alt}_k^0$ holds on $\frk{A}$ iff $\abs{X/E} \leq k$ or, equivalently, $\abs{\frk{A}_0} \leq 2^k$.
\end{proof}

\begin{theorem}
    For any $k \geq 1$ and $\bb{V} \subseteq \bb{MS4} + (\msf{alt}_k^0 = 1)$, $\bb{V}$ is locally finite iff it is of finite depth.
\end{theorem}

\begin{proof}
    We show that in this case \cref{thm:ms4-lf:2} is always satisfied and the criterion reduces to being of finite depth.
    Let $\frk{A} \in \bb{V}$ be s.i.~and $n$-generated. By the above lemma,  $\abs{\frk{A}_0} \leq 2^k$.
    Thus, $\frk{A}_0$ is (at most) $2^k$-generated.
\end{proof}

Since $\bb{S4}_u \models \msf{alt}_1^0$ (by \cref{lem:altk-b0}), this also gives the following generalization of the Segerberg--Maksimova theorem to varieties of $\msf{S4}_u$-algebras:

\begin{corollary}
    \label{cor:s4u-lf-finite-depth}
    $\bb{V} \subseteq \bb{S4}_u$ is locally finite iff it is of finite depth.
\end{corollary}

\begin{remark}
    The formulas $E_i^k$ given in \cite[Thm.~4.2]{Bez02} bound the number of $E_i$-classes in a s.i.~$\msf{S5}^2$-algebra (note that $\bb{S5}^2$ is semisimple, so being s.i.~amounts to being simple) by $k$ or, equivalently, the size of the subalgebra of $\exists_i$-fixpoints by $2^k$.
    Therefore, the local finiteness of the subvarieties of $\bb{S5}^2$ defined by the $E_i^k$ as given in \cite[Lem.~4.4]{Bez02} also follows from our criterion.
\end{remark}

In the general case, however, verifying 
\cref{thm:ms4-lf:2} is more complex. In \cref{translation} we
demonstrate that completely characterizing local finiteness even in $\msfv[2]$ would also do so for $\bb{S5}_2$, which is an open problem as mentioned in the introduction.
We do this by demonstrating a translation from subvarieties of $\bb{S5}_2$ to subvarieties of $\msfv[2]$ that preserves and reflects local finiteness. 
Before that, we will focus our discussion on the behavior of finitely generated $\msf{MS4}$-algebras.

\section{Finitely generated algebras}
\label{finitely-generated-algebras}

In this section we obtain several results about finitely generated $\msf{MS4}$-algebras, including the description of the dual space of the Erd\"os--Tarski algebra. We will see that they are significantly more complex than finitely generated $\msf{S4}$-algebras. We conclude the section with several open problems about finitely generated $\msf{MS4}$-algebras.

We start by recalling the notion of a (quasi-)maximal point (see \cite[Def.~1.4.9]{Esa19}).

\begin{definition}
Let $X$ be a set and $R$ a quasi-order on $X$.
Given $A \subseteq X$, we write 
\[
\max A = \setbuilder{x \in A}{x R y, y \in A \to y R x}
\]
(so for us, $\max$ means ``quasi-maximum'').
We define the \textit{n-th layer} of $X$ inductively as
\[
D_1 = \max X \qquad D_{n+1} = \max \left( X - \bigcup_{i=1}^{n} D_i \right).
\]
\end{definition}

We may speak of the layers of an $\msf{S4}$-frame or $\msf{MS4}$-frame, always referring to the $R$-layers.
In the case of arbitrary quasi-orders, some or all $D_i$ may be empty (e.g., when $X$ contains infinite ascending chains).
However, we have the following (see \cite[Cor.~3.2.2, Thm.~3.2.3]{Esa19}):

\begin{theorem}
Let $(X, R)$ be a descriptive $\msf{S4}$-frame.
Then $D_1 = \max X$ is nonempty and closed.
\end{theorem}

Clearly this holds for descriptive $\msf{MS4}$-frames as well.
We concentrate on finitely generated descriptive frames; that is, frames $\frk{F}$ whose dual $\frk{F}^*$ is a finitely generated algebra. 
Finitely generated descriptive $\msf{S4}$-frames have a relatively well-understood, but complicated structure.
They consist of finite `discrete' layers, possibly with limit points of infinite depth:

\begin{theorem}[see, e.g., {\cite[Sec.~8.6]{CZ97}}]
\label{th:fg-s4}
Let $(X, R)$ be a finitely generated descriptive $\msf{S4}$-frame.
Then every $D_n$ ($n < \omega$) is finite and consists of isolated points, hence is clopen.
\end{theorem}

\begin{remark}
The same result also holds for Heyting algebras and 
their dual spaces; see \cite[Thm.~3.1.8]{Bez-PhD}.
\end{remark}

This result does not carry over to the monadic setting. 
In fact, the Erd\"os--Tarski algebra is an infinite one-generated $\msf{S5}^2$-algebra
(see \cite[Thm.~2.1.11]{HMT85}).
Thus, its dual is an infinite one-generated descriptive $\msf{MS4}$-frame of depth 1. It follows that finitely generated descriptive $\msf{MS4}$-frames may have limit points already in depth 1.
One of the main goals of this section is to describe the dual frame of the Erd\"os--Tarski algebra. 
For this purpose, we 
adjust the notion of a correct partition to the setting of $\msf{MS4}$-frames.

\begin{definition}[\cite{Esa77}]
    \label{def:correct-partition}
    Let $\frk{F} = (X, R, E)$ be an $\msf{MS4}$-frame.
    A \textit{correct partition} of $\frk{F}$ is an equivalence relation $K \subseteq X^2$ such that
    \begin{enumerate}
        \item $K$ is \textit{separated} (or \textit{Boolean}, as in \cite[Def.~8.3]{Kop89}), i.e.~any $x, y$ with $\neg (x K y)$ can be separated by a $K$-saturated clopen set.
        \item $RK \subseteq KR$ ($K$ is \textit{correct} with respect to $R$);
        \item $EK \subseteq KE$ or, equivalently, $EK = KE$ ($K$ is correct with respect to $E$).
    \end{enumerate}
    The quotient $\frk{F}/K$ is defined to be $(X/K, \overline{R}, \overline{E})$, where $\alpha \overline{R} \beta$ iff $\exists x \in \alpha, y \in \beta : x R y$ (or, equivalently, $\forall x \in \alpha, \exists y \in \beta : x R y$), and $\overline{E}$ is defined similarly.
\end{definition}

Correct partitions are the kernels of continuous p-morphisms.
Hence, by standard duality theory (as outlined in \cref{rem:jt-duality}), the lattice of correct partitions of $\frk{F}$ ordered by inclusion is dually isomorphic to the lattice of subalgebras of $\frk{F}^*$. 

\begin{remark}
    \label{rem:correct-partition-equivalence}
    As noted in the previous definition, when a partition $K$ is correct with respect to the equivalence relation $E$, the relations actually commute and we have $EK = KE$.
    Thus, not only does $E$ give rise to an equivalence relation $\overline{E}$ on $K$-classes, but $K$ also gives rise to an 
    equivalence relation on $E$-classes, which we denote by $\overline{K}$.
    In particular, two $E$-classes are $\overline{K}$-related iff every member of each class is $K$-related to some member of the other.
    This idea is put to use in \cref{thm:one-generated} and \cref{lem:partition-lift}.
\end{remark}

We now adapt to our setting a dual characterization of finitely generated algebras; this criterion was given in the setting of $\msf{S4}$-algebras and Heyting algebras in \cite{Esa77}, with an identical proof.

\begin{definition}
    Let $\frk{F} = (X, R, E)$ be an $\msf{MS4}$-frame, and $\overline{g} = (g_1, \dots, g_n)$ a tuple of elements of $\frk{F}^*$.
    Define
    \begin{enumerate}
        \item the coloring $c_{\bar{g}} : X \to 2^n$ by 
        \[
          c_{\bar{g}}(x)(i) = \begin{cases}
            1 & x \in g_i, \\
            0 & \text{otherwise};
          \end{cases}
        \]
        \item the relation $C_{\bar{g}} \subseteq X^2$ as the kernel of $c_{\bar{g}}$, i.e.~$x C_{\bar{g}} y$ iff $c_{\bar{g}}(x) = c_{\bar{g}}(y)$.
    \end{enumerate}
    We call $c_{\bar{g}}$ the \textit{coloring} and $C_{\bar{g}}$ the \textit{partition induced by $\bar{g}$}.
    We will commonly drop the subscript $\bar{g}$ when unambiguous.
\end{definition}

We call an equivalence relation $K$ on $X$ \textit{non-trivial} if it relates two distinct elements of $X$ (i.e., it is strictly above the diagonal in the lattice of equivalence relations).

\begin{theorem}[Coloring Theorem]
\label{thm:coloring-theorem}
Let $\frk{F}$ be a descriptive $\msf{MS4}$-frame, and $\bar{g} = (g_1, \dots, g_n)$ a tuple of elements of $\frk{F}^*$.
Then $\frk{F}^*$ is \textit{not} generated by $\overline{g}$ iff there exists a non-trivial correct partition $K$ of $\frk{F}$ with $K \subseteq C$ (that is, each class of $K$ contains only points of the same color).
\end{theorem}

\begin{proof}
Suppose $\frk{F}^*$ is generated by $\bar{g}$,
and let $K$ be a non-trivial correct partition.
Then $K$ corresponds to a proper subalgebra of $\frk{F}^*$ that couldn't contain all of the $g_i$. 
That is, there is some $g_i$ that is not $K$-saturated.
So there are $x \in g_i$ and $y \not\in g_i$ that are $K$-related. Thus, $xKy$ but $c(x) \neq c(y)$, and we've shown that any non-trivial correct partition $K$ cannot satisfy $K \subseteq C$.

Now suppose $\frk{F}^*$ is not generated by $\bar{g}$.
Then $\bar{g}$ generates a proper subalgebra $\frk{A}$, which gives rise to a non-trivial correct partition $K$, defined by $x K y$ iff for all $a \in \frk{A}$, $x \in a \leftrightarrow y \in a$. 
Now, if $c(x) \neq c(y)$, then $x \in g_i$ and $y \in -g_i$ for some $g_i \in \frk{A}$, so $\neg (xKy)$.
Thus, $K \subseteq C$.
\end{proof}

\begin{remark}
    The previous definition and theorem hold more generally (with the same proof) for any boolean algebra with unary operators
    or, dually, any descriptive frame with binary relations.
\end{remark} 

With all the tools we need, we are now ready to give an explicit description of the $\msf{MS4}$-frame 
that is the dual frame of the Erd\"os--Tarski algebra.

\begin{example}
\label{hmt_compact}
Let $\frk{F} = (X, E_1, E_2)$ where
\[
  X = (\omega + 1 \times \omega + 1) - \set{(\omega, \omega)} \cup \set{k_1, k_2}.
\]
The equivalence classes of $E_1$ are rows, the equivalence classes of $E_2$ are columns, and $k_1, k_2 \in E_1(0, \omega) \cap E_2(\omega, 0)$ (so that $\set{k_1, k_2}$ is the only non-trivial equivalence class of $E_1 \cap E_2$); see \cref{fig:hmt_compact}.
\end{example}

\begin{figure}[h]
    \centering
    \vspace{0.1in}
    \begin{tikzpicture}[
        dot/.style={circle,fill=black,minimum size=6pt,inner sep=0}
    ]
    \foreach \x in {0,...,4}
        \foreach \y in {0,...,4}
            \node [dot] at (\x, \y) {};
    
    \foreach \i in {0,...,4} {
        \node [dot] at (\i, 6) {};
        \node [dot] at (6, \i) {};
    }
    
    \foreach \x in {0,...,4,6} {
        \node at (\x, 5) {\vdots};
        \draw (\x, 0) -- (\x, 4.5);
        \draw (\x, 5.5) -- (\x, 6);
    }
    
    \foreach \y in {0,...,4,6} {
        \node at (5, \y) {\dots};
        \draw (0, \y) -- (4.5, \y);
        \draw (5.5, \y) -- (6, \y);
    }
    
    \draw[rounded corners,stroke=black,fill=white] (5.6, 5.6) rectangle (6.4, 6.4);
    
    \node [dot,label=west:$k_1$] at (5.8, 6.2) {};
    \node [dot,label=east:$k_2$] at (6.2, 5.8) {};
    
    \draw[rounded corners,stroke=black,opacity=0.5] (-0.35, -0.55) -- (-0.35, 6.45) -- (6.65,6.45) -- cycle;
    
    \node at (-0.55, 6) {$g$};
    \end{tikzpicture}
    \caption{Descriptive frame of \cref{hmt_compact}.}
    \label{fig:hmt_compact}
    \end{figure}
    
In the remainder of the section, we demonstrate that $\frk{F}$ has the desired properties.
To define the topology on $X$, we will describe it as a two-point compactification of the space $X_0 = (\omega + 1 \times \omega + 1) - \set{(\omega, \omega)}$, where $\omega + 1$ has the usual interval topology, and $X_0$ has its usual topology of a subspace of the product.
As an open subset of a 
Stone space, $X_0$ is zero-dimensional, locally compact, and Hausdorff (that is, $X_0$ is locally Stone). 
Define
\[
  g_1 = \setbuilder{(i, j)}{i \leq j \leq \omega} \cap X_0 \qquad g_2 = X_0-g_1.
\]
Note that $g_1$ and $g_2$ are disjoint clopen sets that are neither compact nor co-compact (i.e., having a compact complement) in $X_0$.
In the terminology of \cite[Def.~2.2]{Mag65}, $\set{g_1, g_2}$ is a \textit{2-star} in $X_0$---a pair of disjoint open sets with the property that the complement of their union is compact but the complement of each set itself is not.
In \cite[Thm.~2.1]{Mag65} it is shown that this datum gives rise to a 2-point compactification of $X_0$; specifically, we have the following basis for a compact Hausdorff topology on $X$:
\[
  \mca{B} = \Clp(X_0) \cup \setbuilder{U \cup \set{k_i}}{U \in \Clp(X_0) \tand g_i - U \text{ is compact in $X_0$}}
\]
(Note that $\mca{B}$ is not a clopen basis.)

\begin{lemma}
    $X$ is a Stone space.
\end{lemma}

\begin{proof}
    By \cite[Thm.~2.1]{Mag65}, 
    $X$ is compact Hausdorff. 
    We show that 
    $X$ is zero-dimensional. 
    For this it suffices to separate distinct points of $X$ by a clopen subset of $X$ (see, e.g., \cite[69]{Joh82}). 
    First suppose that $x, y \in X_0$. Since $X_0$ is zero-dimensional, there is a clopen set $U\subseteq X_0$ containing $x$ and missing $y$. Because $X_0$ is locally compact, $x$ has a compact neighborhood $K$ in $X_0$. Therefore, there is a clopen set $V$ with $x\in V \subseteq U\cap K$. But then $V$ is 
    a compact clopen of $X_0$, and hence a clopen subset of $X$ separating $x$ and $y$.
    The points $k_1$ and $k_2$ are separated by $g_1 \cup \set{k_1}$ 
    which is clopen of $X$ by definition.
    Finally, $x \in X_0$ and $k_i$ are separated by any compact clopen neighborhood $K$ of $x$ in $X_0$, which is a clopen set of $X$. 
\end{proof}

\begin{lemma}
    $\frk{F}$ is a descriptive $\msf{S5}^2$-frame.
\end{lemma}

\begin{proof}
    Since $X$ is a Stone space by the above lemma,
    it remains to show that $E_1, E_2$ are 
    continuous. We show that $E_1$ is continuous. The proof for $E_2$ is symmetric.
    Observe that all finite-index rows (i.e., $[0, \omega] \times \set{j}$ for $j < \omega$) are evidently clopen, while the ``top row'' $\setbuilder{(i, \omega)}{i < \omega} \cup \set{k_1,k_2}$ is closed but not open (its complement is the union of all finite-index rows).
    Hence, $E_1$ is point-closed. 

    It is left to show that $U$ clopen implies $E_1(U)$ is clopen. By \cite[Thm.~3.1.2(IV)]{Esa19}, it is enough to verify the following two conditions:
    \begin{enumerate}
        \item If $y \not\in E_1(x)$, then there is a partition into two $E_1$-saturated open sets separating $x$ and $y$.
        \item For all $x$ and open $V$, if $E_1(x) \cap V \neq \varnothing$, then there is an open neighborhood $U$ of $x$ with $E_1(y) \cap V \neq \varnothing$ for all $y \in U$.
    \end{enumerate}
    For (1), since $y \not\in E_1(x)$, $x$ and $y$ are in distinct rows. Therefore, one of them must be in a row of finite index, say $R_n = [0,\omega] \times \{ n \}$ with $n < \omega$ (if both are in rows of finite index, take $n$ to be the least index).
    Then $[0, \omega] \times \set{n}$ and its complement provide a partition of $X$ into two $E_1$-saturated open sets separating $x$ and $y$.
    
    For (2), we consider cases:
    If $x = (i, j)$ for $i, j < \omega$, $\set{x}$ is an appropriate neighborhood.
    If $x = (\omega, j)$ for $j < \omega$,
    then $[0,\omega] \times \set{j}$ is an appropriate neighborhood.
    Suppose $x = (i, \omega)$ for $i<\omega$.
    If $E_1(x) \cap V$ contains $(i', \omega)$ for $i'<\omega$, then $V$ contains a neighborhood $\set{i'} \times [n, \omega]$ for $n < \omega$, so $\set{i} \times [n, \omega]$ is an appropriate neighborhood of $x$.
    If $E_1(x) \cap V$ contains $k_1$, then $V$ contains some $V' \cup \set{k_1}$ with $g_1 - V'$ compact in $X_0$, as sets of this form are the only basic open neighborhoods of $k_1$.
    Then there must be some $n$ with $(g_1- V')\cap [n, \omega]^2 = \varnothing$ (indeed any $U \subseteq X_0$ compact 
    has this property since otherwise $\setbuilder{X_0-[n, \omega]^2}{n < \omega}$ would be an open cover of $U$ without a finite subcover).
    Therefore, $g_1 \cap [n, \omega]^2 \subseteq V' \subseteq V$,
    and again $\set{i} \times [n, \omega]$ is an appropriate neighborhood of $x$.
    If $E_1(x) \cap V$ contains $k_2$, then by the same argument we have an $n$ with $g_2 \cap [n, \omega]^2 \subseteq V$ and the same neighborhood works.
    Finally, suppose $x = k_s$ for $s \in \set{1,2}$.
    If $E_1(x) \cap V$ contains $(i, \omega)$ for $i < \omega$, then $V$ contains $\set{i} \times [n, \omega]$ for some $n$, and $(g_s \cap [n, \omega]^2) \cup \set{k_s}$ is an appropriate basic open neighborhood of $x$.
    If $E_1(x) \cap V$ contains $k_t$ for $t \neq s$, then by the same argument as before we have $g_t \cap [n, \omega]^2 \subseteq V$ for some $n$, and the same neighborhood $(g_s \cap [n, \omega]^2) \cup \set{k_s}$ works.
\end{proof}

\begin{theorem}
\label{thm:one-generated}
$\frk{F}$ is one-generated and $\frk{F}^*$ is isomorphic to the Erd\"os--Tarski algebra.
\end{theorem}

\begin{proof}
Let $g = g_1 \cup \set{k_1}$ (see \cref{fig:hmt_compact}).
We claim $\frk{F}^* = \la g \ra$.
Suppose not. By the coloring theorem, there exists a non-trivial correct partition $K$ of $X$ such that 
every class of $K$ is contained in $g$ or $-g$.
Since it is non-trivial, $K$ must identify two points in distinct rows or columns (by assumption we cannot have $k_1 K k_2$).
By \cref{rem:correct-partition-equivalence}, this is to say that two distinct rows or columns must be $\overline{K}$-related.
Let $C_i = \set{i} \times [0, \omega]$ be the $i$-th column and 
$R_j = [0, \omega] \times \set{j}$ the $j$-th row.
Observe that $C_0$
cannot be related to $C_i$
for $i > 0$ 
since otherwise a point in $-g \cap C_i$ would be $K$-related to a point in $-g \cap C_0 = \varnothing$.

Now suppose $C_{i_1}$ and $C_{i_2}$
are $\overline{K}$-related for $0 < i_1 < i_2 \leq \omega$.
Then there must be some $j_2 \in [i_1, i_2)$ so that $(i_2, j_2) K (i_1, j_1)$ for some $j_1 \in [0, i_1)$ (since all points in $-g \cap C_{i_2}$ must be $K$-related to points in $-g \cap C_{i_1}$).
Therefore, rows $j_1$ and $j_2$ are $\overline{K}$-related, with $0 \leq j_1 < j_2 < i_2$.
Thus, $(j_2, j_2) K (i_0, j_1)$ for some $i_0 \in [0, j_1]$.
We now conclude that columns $i_0$ and $j_2$ are $\overline{K}$-related, and note that $0 \leq i_0 < i_1$, $0 < j_2 < i_2$, and $i_0 < j_2$.
So from our assumption of two $\overline{K}$-related columns, we have found two $\overline{K}$-related columns of strictly lower index, and iterating this process will eventually reveal that column $0$ must be $\overline{K}$-related to some non-zero column, which is a contradiction by the previous paragraph.

Finally, suppose $R_{j_1}$ and $R_{j_2}$
are $\overline{K}$-related for $0 \leq j_1 < j_2 \leq \omega$.
Then $(j_2, j_2) K (i, j_1)$ (or $k_1 K (i, j_1)$ in the case $j_2 = \omega$) for some $i \in [0, j_1]$.
But then columns $i$ and $j_2$ are $\overline{K}$-related, with $i < j_2$, which is a contradiction by the last paragraph.

Consequently, such a $K$ cannot exist, 
and we conclude that $\frk{F}^*$ is generated by $g$.
There is an evident homomorphism $f : \frk{F}^* \to \mca{P}(\omega \times \omega)$ given by $U \mapsto U \cap (\omega \times \omega)$.
    Then $f(g) = \setbuilder{(i, j)}{i \leq j}$, and $f$ is an isomorphism onto its image (because $\omega\times\omega$ is dense in $X$).
    Thus, $\mathfrak{F}^*$ is isomorphic to the subalgebra of $\mathcal{P}(\omega \times \omega)$ generated by $f(g)$, which is exactly the Erd\"os--Tarski algebra.
\end{proof}

This provides an example of a descriptive $\msf{MS4}[1]$-frame whose (only) layer is infinite and contains limit points.
Thus, the only part of \cref{th:fg-s4} that has a hope to be salvaged in the monadic setting 
is that each layer may be clopen:

\begin{question}
Let $(X, R, E)$ be a finitely generated descriptive $\msf{MS4}$-frame.
Is $D_1 = \max X$ or any $D_n$ a clopen (admissible, definable) set?
\end{question}

\begin{question}
Does there exist a finitely generated descriptive $\msf{MS4}$-frame whose finite layers consist \textit{only} of limit points?
\end{question}

\section{Translating $\bb{S5}_2$ to $\msfv[2]$}
\label{translation}

In this section we follow up on the promise made at the end of \cref{local-finiteness} and give a semantic translation from subvarieties of $\bb{S5}_2$ to subvarieties of $\msfv[2]$ that preserves and reflects local finiteness.
Before doing so, we develop some basic facts about the variety of $\msf{S5}_2$-algebras and their dual frames.

\begin{definition}[{\cite[Sec.~3.1]{GKWZ03}}]
    The logic $\msf{S5}_2$ is the {\em fusion} $\msf{S5} \oplus \msf{S5}$, i.e.~the smallest normal modal
    logic in the language with two modalities $\exists_1$ and $\exists_2$ containing the $\msf{S5}$ axioms for each $\exists_i$ (and no other axioms). 
\end{definition}

Algebraic models of $\msf{S5}_2$ are triples $(B, \exists_1, \exists_2)$, where $(B, \exists_i)$ is an $\msf{S5}$-algebra for $i = 1,2$. We call such algebras {\em $\msf{S5}_2$-algebras}. 
J\'{o}nsson--Tarski duality specializes to yield that the corresponding variety $\bb{S5}_2$ is dually equivalent to the following category of descriptive frames:

\begin{definition} 
    A {\em descriptive $\msf{S5}_2$-frame} is a triple $\frk{F} = (X, E_1, E_2)$ where $X$ is a Stone space and $E_1,E_2$ are two continuous equivalence relations on $X$ (see \cref{def: cont rel}).
\end{definition}

An {\em $\msf{S5}_2$-morphism} between descriptive $\msf{S5}_2$-frames is a continuous map $f:(X, E_1, E_2) \to (X', E_1', E_2')$ that is a p-morphism with respect to both $(E_1,E_1')$ and $(E_2,E_2')$.
When we depict $\msf{S5}_2$-frames, we will depict $E_1$-equivalence classes as horizontal lines and $E_2$-equivalence classes with vertical analogues in blue.

\begin{definition}\
    \begin{enumerate}
        \item For an $\msf{S5}_2$-algebra $(B, \exists_1, \exists_2)$, define $\diamondsuit = \exists_1 \vee \exists_2$.
        \item For an $\msf{S5}_2$-frame $(X, E_1, E_2)$, define $S = E_1 \cup E_2$.
    \end{enumerate}
\end{definition}

It is straightforward to see that $S$ is a reflexive and symmetric relation (or \textit{edge relation}) on $X$ and that $S$ is the dual relation of $\diamondsuit$.
Hence, $\diamondsuit$ is a possibility operator validating the well-known $\msf{T}$ and $\msf{B}$ axioms, making it a $\msf{KTB}$-modality. 
We write $S^* = \bigcup_{n=0}^\infty S^n$ for the reflexive transitive closure of $S$.

\begin{definition}[\cite{Ven04}]
    Let $\frk{F} = (X, E_1, E_2)$ be a descriptive $\msf{S5}_2$-frame.
    \begin{enumerate}
    \item Call $x \in X$ a \textit{topo-root} of $\frk{F}$ if $S^*(x)$ is dense in $X$. 
    \item We say that $\frk{F}$ is \textit{topo-rooted} if it has a topo-root, and \textit{strongly topo-rooted} if it has an open set of topo-roots.
    \end{enumerate}
\end{definition}

The following is a consequence of general results in \cite{Ven04}:

\begin{theorem} 
    Let $\frk{A}$ be an $\msf{S5}_2$-algebra.
    \begin{enumerate}
        \item $\frk{A}$ is 
        s.i.~iff $\frk{A}_*$ is strongly topo-rooted.
        \item $\frk{A}$ is simple iff every point of $\frk{A}_*$ is a topo-root.
    \end{enumerate}
\end{theorem}

Following \cite{Rau80} (see also \cite[74]{Kra99}), we call a variety $\bb{V} \subseteq \bb{S5}_2$ \textit{$n$-transitive} if $\bb{V} \models \diamondsuit^{n+1} a \leq \diamondsuit^n a$, and \textit{weakly transitive} if it is $n$-transitive for some $n < \omega$.
The following is a consequence of a more general result of Kowalski and Kracht \cite[Thm.~12]{Kow06} (note that $\bb{S5}_2$ is automatically \textit{cyclic} since each basic modality is $\msf{S5}$, see \cite[Prop.~6]{Kow06}).

\begin{theorem}
    A variety $\bb{V} \subseteq \bb{S5}_2$ is semisimple iff it is weakly transitive.
\end{theorem}

It is a consequence of \cite[Cor.~1]{Ven04} (where weak transitivity goes by the name of \textit{$\omega$-transitivity}) that in a weakly transitive variety of $\msf{S5}_2$-algebras, roots and topo-roots coincide (where $x$ is a root of $\mathfrak{F}$ if $S^*(x) = X$); together with the previous theorem, this yields:

\begin{theorem} 
    Let $\bb{V} \subseteq \bb{S5}_2$ be a weakly transitive variety. The following are equivalent:
    \begin{enumerate}
        \item $\frk{A} \in \bb{V}$ is s.i.
        \item $\frk{A} \in \bb{V}$ is simple.
        \item Every point of $\frk{A}_*$ is a root.
    \end{enumerate}
\end{theorem}

We are now ready to describe our translation. We start by 
a construction producing a descriptive $\msfl[2]$-frame from a descriptive $\msf{S5}_2$-frame.

\begin{construction}
Let $\frk{F} = (X, E_1, E_2)$ be a descriptive $\msf{S5}_2$-frame. We let $\mca{L} = X / E_2$ be the quotient space and $\pi_\mca{L} : X \to \mca{L}$ the quotient map.
We set $T(\frk{F}) = (X \cup \mca{L}, R, E)$, where
\begin{itemize}
    \item $X \cup \mca{L}$ is the (disjoint) union of $X$ and $\mca{L}$,
    \item $E$ is the smallest equivalence on $X \cup \mca{L}$ containing $E_2$ along with
\[
\setbuilder{(x, \alpha)}{x \in X, \alpha \in \mca{L}, x \in \alpha},
\]
    \item $R = E_1 \cup (X \times \mca{L}) \cup (\mca{L} \times \mca{L})$.
\end{itemize}
In the case that $\abs{\mathfrak{F}}$ is finite, it follows that $\lvert T(\mathfrak{F}) \rvert \leq 2\lvert \mathfrak{F} \rvert$ and $\lvert T(\mathfrak{F})^* \rvert \leq \lvert \mathfrak{F}^* \rvert^2$.
\end{construction}

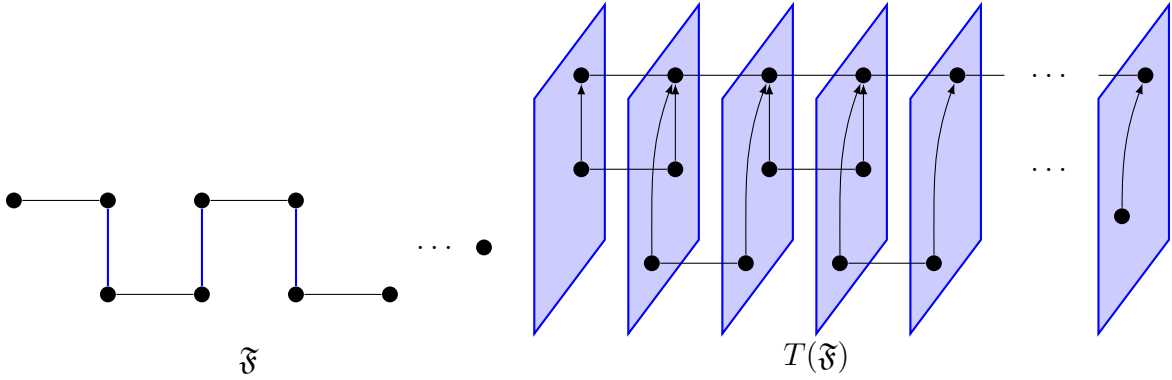
\begin{figure}[h]
\centering
\begin{tikzpicture}[
    scale=1.25,
	dot/.style={circle,fill=black,minimum size=6pt,inner sep=0}
]
\node (A) [dot] at (0, 1) {};
\node (B) [dot] at (1, 1) {};
\node (C) [dot] at (1, 0) {};
\node (D) [dot] at (2, 0) {};
\node (E) [dot] at (2, 1) {};
\node (F) [dot] at (3, 1) {};
\node (G) [dot] at (3, 0) {};
\node (H) [dot] at (4, 0) {};

\draw (A) -- (B);
\draw[blue,thick] (B) -- (C);
\draw (C) -- (D);
\draw[blue,thick] (D) -- (E);
\draw (E) -- (F);
\draw[blue,thick] (F) -- (G);
\draw (G) -- (H);

\node at (4.5, 0.5) {$\dots$};
\node [dot] at (5, 0.5) {};

\node  at (2.5, -0.7) {$\frk{F}$};
\end{tikzpicture}
\quad
\begin{tikzpicture}[
    scale=1.25,
	dot/.style={circle,fill=black,minimum size=6pt,inner sep=0}
]

\foreach \x in {0,...,4,6} {
	\draw[blue,thick,fill=blue!20] (\x - 0.5, -0.75) -- (\x - 0.5, 1.75) -- (\x + 0.25, 2.75) -- (\x + 0.25, 0.25) --cycle;
}

\node (A2) [dot] at (0, 1) {};
\node (B2) [dot] at (1, 1) {};
\node (C2) [dot] at (0.75, 0) {};
\node (D2) [dot] at (1.75, 0) {};
\node (E2) [dot] at (2, 1) {};
\node (F2) [dot] at (3, 1) {};
\node (G2) [dot] at (2.75, 0) {};
\node (H2) [dot] at (3.75, 0) {};

\node (I2) [dot] at (5.75, 0.5) {};

\node (A1) [dot] at (0, 2) {};
\node (B1) [dot] at (1, 2) {};
\node (C1) [dot] at (2, 2) {};
\node (D1) [dot] at (3, 2) {};
\node (E1) [dot] at (4, 2) {};

\node (F1) [dot] at (6, 2) {};

\draw (A2) -- (B2);
\draw (C2) -- (D2);
\draw (E2) -- (F2);
\draw (G2) -- (H2);

\draw (A1) -- (E1) -- (4.5, 2);
\draw (5.5, 2) -- (F1);
\node at (5, 2) {$\dots$};

\draw[-latex] (A2) -- (A1);
\draw[-latex] (B2) -- (B1);
\draw[-latex] (E2) -- (C1);
\draw[-latex] (F2) -- (D1);

\draw[-latex] (C2) to[out=90, in=250] (B1);
\draw[-latex] (D2) to[out=90, in=250] (C1);
\draw[-latex] (G2) to[out=90, in=250] (D1);
\draw[-latex] (H2) to[out=90, in=250] (E1);

\draw[-latex] (I2) to[out=90, in=250] (F1);

\node at (5, 1) {$\dots$};

\node at (2.5, -1) {$T(\frk{F})$};
\end{tikzpicture}
\caption{
    Constructing an $\msfl[2]$-frame from an $\msf{S5}_2$-frame.
    In the $\msf{S5}_2$-frame on the left, $E_1$-clusters are horizontal lines while $E_2$-clusters are blue vertical lines.
    In the $\msfl[2]$-frame on the right, $R$-clusters are horizontal lines, proper $R$-arrows are drawn with arrowheads, and $E$-clusters are given by the blue rectangles.
}
\label{fig:translation}
\end{figure}

This construction is depicted in \cref{fig:translation}.
Under the definition of $R$, $\mca{L}$ becomes an $R$-cluster, which we think of as a ``top rail''.

To prove that $T(\frk{F})$ is a descriptive $\msfl[2]$-frame, we need the following:

\begin{lemma} \label{lem: clopen map}
    $\pi_\mca{L}$ is a clopen map, meaning that
    $U$ clopen implies $\pi_\mca{L}[U]$ is clopen. 
\end{lemma}

\begin{proof}
    Suppose $U \subseteq X$ is clopen.
    Then $\pi_\mca{L}^{-1}(\pi_\mca{L}[U]) = E_2(U)$, which is clopen by continuity of $E_2$.
    Since $\pi_\mca{L}$ is a quotient map, $\pi_\mca{L}[U]$ is clopen.
\end{proof}

\begin{lemma}
    \label{lem:translation-descriptive-frame}
    $T(\frk{F})$ is a descriptive $\msfl[2]$-frame; furthermore, $T(\frk{F})$ is $Q$-rooted.
\end{lemma}

\begin{proof}
    Note that $E_2$ is 
    a separated partition of $X$ (but in general not correct with respect to $E_1$ in the sense of \cref{def:correct-partition}) since any $(x, y) \not\in E_2$ can be separated by a $E_2$-saturated clopen (this is a property of any continuous quasi-order; see \cite[Thm.~3.1.2]{Esa19}).
    This ensures $\mca{L}$ is a Stone space (\cite[Lem.~8.4]{Kop89}). Thus, so is $X \cup \mca{L}$ (the union is disjoint).

    Clopen sets in $T(\frk{F})$ are of the form $U \cup V$, with $U$ clopen in $X$ and $V$ clopen in $\mca{L}$.
    We have
    \begin{itemize}
        \item $E(x) = E_2(x) \cup \set{E_2(x)}$ for $x \in X$;
        \item $E(\alpha) = \alpha \cup \set{\alpha}$ for $\alpha \in \mca{L}$; 
        \item $E(U \cup V) = E_2(U) \cup \pi_\mca{L}[U] \cup V \cup \pi_\mca{L}^{-1}(V)$.
    \end{itemize}
    The first two items show that $E$-clusters are 
    closed.
    In the third item, every set on the right hand side is clopen, by continuity of $E_2$ and the fact that $\pi_\mca{L}$ is a clopen map (\cref{lem: clopen map}). Thus, $E$ is continuous.
    For $R$, we have 
    \begin{itemize}
        \item $R(x) = E_1(x) \cup \mca{L}$ for $x \in X$;
        \item $R(\alpha) = \mca{L}$ for $\alpha \in \mca{L}$;
    \item $
    R^{-1}(U \cup V) = \begin{cases}
        E_1(U) & V = \varnothing, \\
        X \cup \mca{L} & \text{otherwise}.
    \end{cases}
    $
    \end{itemize}
    By the first two items, $R$ is point-closed, which together with the third yields that $R$ is continuous.

    Now observe that $Q = ER$ is the total relation on $T(\frk{F})$ (i.e., $Q = (X \cup \mca{L})^2$) because for each $x\in X$ we have $\mca{L}\subseteq R(x)$, so $X \cup \mca{L} = ER(x)$.
    Therefore, every point of $T(\frk{F})$ is a $Q$-root and, in particular, 
    $RE \subseteq ER$.
    Finally, $T(\frk{F})$ is of $R$-depth 2 by construction. Thus, $T(\frk{F})$ is a $Q$-rooted descriptive $\msfl[2]$-frame.
\end{proof}

By construction, $T(\frk{F})$ has exactly two layers $D_1 = \mca{L}$ and $D_2 = X$.
By design, we have $E \vert_{D_2} = E_2$ and $R \vert_{D_2} = E_1$.
Since $D_1, D_2$ are clopen, 
we can consider the \textit{relativization}
\[
\frk{A}_i = (\Clp D_i, \cap, \cup, (-)', \varnothing, D_i, \lozenge', \exists'),
\]
where the relative complementation and modal operators are given by
\[
a' = D_i - a \qquad \lozenge' a = \lozenge a \cap D_i \qquad \exists' a = \exists a \cap D_i.
\]
The maps $a \mapsto a \cap D_i$ are Boolean homomorphisms.
Considering the restrictions of $R$ and $E$ to each layer, we see that $\frk{A}_1$ is an $\msf{S5}^2$-algebra (and in fact belongs to a proper and hence locally finite subvariety of $\bb{S5}^2$), while $\frk{A}_2$ is an $\msf{S5}_2$-algebra isomorphic to $\frk{F}^*$:

\begin{lemma}
    \label{lem:second-layer}
    $\frk{A}_2 \cong \frk{F}^*$, and the relative operations of $\frk{A}_2$ are definable in $T(\frk{F})^*$. 
\end{lemma}

\begin{proof}
    This follows from the explicit description of $R^{-1}$ and $E$ on clopen sets in \cref{lem:translation-descriptive-frame}.
    Also, the relative complementation and modal operators are definable using the clopen set $D_2 \in T(\frk{F})^*$, e.g.~by using precisely the definitions from the preceding paragraph.
\end{proof}

Recall from \cref{rem:correct-partition-equivalence} that, since $E_2$ is an equivalence, a correct partition $K$ on $\frk{F}$ yields an equivalence relation $\overline{K}$ on the set $\mca{L}$ of $E_2$-classes of $\frk{F}$.

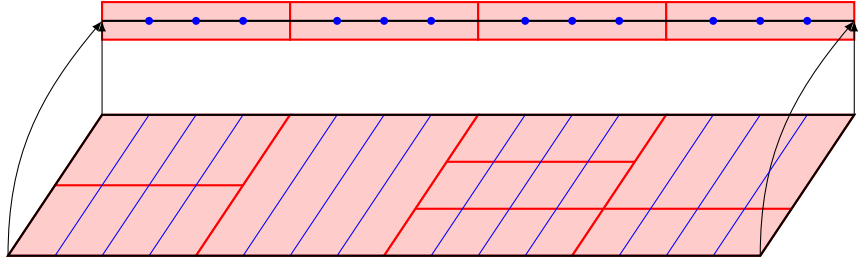
\begin{figure}[h]
\centering
\begin{tikzpicture}[
    scale=2.5,
	dot/.style={circle,fill=black,minimum size=3pt,inner sep=0}
]
\foreach \x in {0,...,3} {
	\draw[red,thick,fill=red!20] (\x - 0.5, -0.75) -- (\x, 0) -- (\x + 1, 0) -- (\x + 0.5, -0.75) --cycle;
}
\foreach \y in {-0.75/2} {
    \draw[red,thick] (0 + 2*\y/3, \y) -- (0 + 2*\y/3 + 1, \y);
}
\foreach \y in {-0.75/3, -2*0.75/3} {
    \draw[red,thick] (2 + 2*\y/3, \y) -- (2 + 2*\y/3 + 1, \y);
}
\foreach \y in {-2*0.75/3} {
    \draw[red,thick] (3 + 2*\y/3, \y) -- (3 + 2*\y/3 + 1, \y);
}
\foreach \x in {0,...,3} {
	\draw[red,thick,fill=red!20] (\x, 0.4) -- (\x, 0.6) -- (\x + 1, 0.6) -- (\x + 1, 0.4) --cycle;
}

\draw[thick] (0, 0.5) -- (4, 0.5);

\draw[-latex] (0, 0) -- (0, 0.5);
\draw[-latex] (-0.5, -0.75) to[out=90,in=230] (0, 0.5);
\draw[-latex] (4, 0) -- (4, 0.5);
\draw[-latex] (3.5, -0.75) to[out=90,in=230] (4, 0.5);

\foreach \x in {0,...,3} {
    \draw[blue] (\x + 0.25, 0) -- (\x - 0.25, -0.75);
    \draw[blue] (\x + 0.5, 0) -- (\x, -0.75);
    \draw[blue] (\x + 0.75, 0) -- (\x + 0.25, -0.75);
    \node [dot,fill=blue] at (\x + 0.25, 0.5) {};
    \node [dot,fill=blue] at (\x + 0.5, 0.5) {};
    \node [dot,fill=blue] at (\x + 0.75, 0.5) {};
}

\draw[thick] (0, 0) -- (-0.5, -0.75) -- (3.5, -0.75) -- (4, 0) -- cycle;

\end{tikzpicture}
\caption{
    Lifting a correct partition of $\frk{F}$ to one of $T(\frk{F})$, as in \cref{lem:partition-lift}:
    $E$-classes are depicted in blue as $E_2$-classes of $\frk{F}$ with the $E_2$-class itself above on the top rail. 
    The correct partition $\what{K}$ is depicted in red.
}
\label{fig:partition-lift}
\end{figure}

\begin{lemma}
    \label{lem:partition-lift}
    Let $K$ be a correct partition of $\frk{F}$.
    Define $\what{K} = K \cup \overline{K}$, a relation on $T(\frk{F})$.
    That is, $\what{K}$ is the relation $K$ on $D_2 = X$, and the relation $\overline{K}$ on $D_1 = \mca{L}$ (and relates no points on different layers); see \cref{fig:partition-lift}.
    \begin{enumerate}
        \item $\what{K}$ is a correct partition of $T(\frk{F})$.
        \item $T(\frk{F}/K) \cong T(\frk{F})/\what{K}$.
    \end{enumerate}
\end{lemma}

\begin{proof}
    (1) We first show $\what{K}$ is separated.
    If $x, y \in D_2$ and $\neg (x \what{K} y)$, then $\neg (x K y)$ and they can be separated by a $K$-saturated clopen subset of $D_2$, which is also $\what{K}$-saturated (since $\what{K}$ relates no points in different layers).
    Any $x \in D_2$ and $\alpha \in D_1$ can be separated by $D_2$
    which is $\what{K}$-saturated by the same reason.
    Now suppose $\alpha, \beta \in D_1$ and $\neg (\alpha \what{K} \beta)$.
    Then $\neg (x K y)$ for any $x \in \alpha, y \in \beta$, and so any such $x, y$ can be separated by a $K$-saturated clopen $A\subseteq D_2$.
    By \cref{lem: clopen map}, $\pi_\mca{L}[A]$ 
    is a clopen subset of $D_1$ separating $\alpha$ and $\beta$.
    To see it is
    $\what{K}$-saturated, let $\gamma \in \pi_\mca{L}[A]$ and $\gamma \what{K} \delta$.
    Then there is $x \in A \cap \gamma$ and $y \in \delta$ with $x K y$ ($x$ may be chosen from $A$ by the remarks at the end of \cref{def:correct-partition}). Since $A$ is $K$-saturated, $y \in A \cap \delta$, and hence $\delta \in \pi_\mca{L}[A]$.

    We now show $\what{K}$ is correct with respect to $R$.
    Suppose $x, y \in D_2$. 
    If $x \what{K} y R z$ for $z \in D_2$, then 
    $x K y E_1 z$, so $x E_1 y' K z$ by correctness of $K$, and hence $x R y' \what{K} z$.
    If $x \what{K} y R \alpha$ for $\alpha \in D_1$, then $x R \alpha \what{K} \alpha$ (by definition of $R$).
    Suppose $\alpha, \beta \in D_1$.
    If $\alpha \what{K} \beta R \gamma$ for $\gamma \in D_1$ (the only possibility), then $\alpha R \gamma \what{K} \gamma$ (again by definition of $R$).
    In any case, $R\what{K} \subseteq \what{K} R$.

    Finally, we show $\what{K}$ is correct with respect to $E$.
    Suppose $x, y \in D_2$.
    If $x \what{K} y E z$ for $z \in D_2$, then 
    $x K y E_2 z$, so $x E_2 y' K z$, and hence $x E y' \what{K} z$.
    If $x \what{K} y E \alpha$ for $\alpha \in D_1$, then $\alpha = E_2(y)$, and $E_2(x) \, \overline{K} \, E_2(y)$.
    So $x \; E \; E_2(x) \; \what{K} \; \alpha$.
    Suppose $\alpha, \beta \in D_1$.
    If $\alpha \what{K} \beta E z$ for $z \in D_2$ (the only non-trivial possibility), then $\beta = E_2(z)$.
    So we must have $z K x$ for some $x \in \alpha$, and hence $\alpha E x \what{K} z$.
    In any case, $E\what{K} \subseteq \what{K} E$. We conclude that $\what{K}$ is a correct partition of $T(\frk{F})$.

    (2) We exhibit an isomorphism $f : T(\frk{F}/K) \to T(\frk{F})/\what{K}$.
    Let $f(K(x)) = K(x)$ on $K$-classes (i.e., the identity on the common second layer of both frames) and $f(\overline{E_2} (K(x))) = \overline{K} (E_2(x))$ on the additional points.
    To check continuity, it suffices to check on clopen subsets of the top layer only (since $\what{K}$ identifies no points across layers and $f$ is the identity on the bottom layers); such a clopen can be identified with a $\overline{K}$-saturated clopen $U$ of $\frk{F}$, and we have
    \begin{align*}
        f^{-1}(U) =& \setbuilder{\overline{E_2}(K(x))}{E_2(x) \in U} \text{ is clopen in $T(\frk{F}/K)$} \\
        \text{iff } & \setbuilder{K(x)}{E_2(x) \in U} \text{ is clopen in $\frk{F}/K$} \\
        \text{iff } & \setbuilder{x}{E_2(x) \in U} \text{ is clopen in $\frk{F}$}.
    \end{align*}
    The last set is precisely $\pi_\mca{L}^{-1}(U)$, which is clopen since $U$ is.
    Clearly $f$ is surjective; for injectivity, suppose $\overline{E_2}(K(x)) \neq \overline{E_2}(K(y))$.
    Then
    \[
    \neg (K(x) \, \overline{E_2} \, K(y)) \rightarrow \neg (x K y) \rightarrow \neg (E_2(x) \, \overline{K} \, E_2(y)),
    \]
    and hence $\overline{K}(E_2(x)) \neq \overline{K}(E_2(y))$.
    
    Since $f$ is a continuous bijection, it remains to show it preserves and reflects the relations (\cite[Prop.~1.4.15]{Esa19}).
    For $R$, note (checking only the cases where $f$ is not the identity on some input):
    \begin{itemize}
        \item $K(x) \; R \; E_2(K(y)) \leftrightarrow \what{K}(x) \; R \; \what{K}(E_2(y))$ (both are always true by definition).
        \item $E_2(K(x)) \; R \; E_2(K(y)) \leftrightarrow \what{K}(E_2(x)) \; R \; \what{K}(E_2(y))$ (similarly).
    \end{itemize}
    For $E$ we have
    \[
    K(x) \; E \; \overline{E_2}(K(y)) \leftrightarrow K(x) \, \overline{E_2} \, K(y) \leftrightarrow x \, E_2 \, y \leftrightarrow x \; E \; E_2(y) \leftrightarrow \what{K}(x) \; E \; \what{K}(E_2(y)).
    \]
    Thus, $f$ is an isomorphism of $\msf{MS4}$-frames.
\end{proof}

\begin{lemma}
    \label{lem:loc-fin-right-direction}
    Let $g_1, \dots, g_n \in T(\frk{F})^*$ and $L$ be the correct partition of $T(\frk{F})$ corresponding to the $n$-generated subalgebra $\frk{B} = \la g_1, \dots, g_n \ra_{T(\frk{F})^*}$ of $T(\frk{F})^*$.
    Let $K$ be the correct partition of $\frk{F}$ corresponding to the $2n$-generated subalgebra 
    \[
    \frk{B}' = \la g_1 \cap D_2, E(g_1 \cap D_1) \cap D_2, \dots, g_n \cap D_2, E(g_n \cap D_1) \cap D_2 \ra_{\frk{F}^*}
    \]
    of $\frk{F}^*$. Then $\what{K} \subseteq L$.
\end{lemma}

\begin{proof}
    It is sufficient to show that each generator of $\frk{B}$ is $\what{K}$-saturated.
    Since $\what{K}$ relates no points across layers, and each $g_i \cap D_2$ is $K$-saturated (as a generator of $\frk{B}'$), we need only check that $g_i \cap D_1$ is $\overline{K}$-saturated.
    Suppose $\alpha \in g_i$ and $\alpha \overline{K} \beta$.
    Then there are $x \in \alpha, y \in \beta$ with $x K y$.
    Evidently, $x \in E(g_i \cap D_1) \cap D_2$ which is $K$-saturated (as a generator of $\frk{B}'$), so $y \in E(g_i \cap D_1) \cap D_2$.
    That is, $y$ is $E$-related to something in $g_i \cap D_1$, which could only be $\beta$, and hence $\beta \in g_i$.
\end{proof}

We now adapt our translation on frames to one on varieties of algebras.
For $\frk{A} \in \bb{S5}_2$, let
\[
  T(\frk{A}) = (T(\frk{A}_*))^*.
\]
We write $H$, $S$, and $P$ for the class operators of taking homomorphic images, subalgebras, and products, respectively.
For a variety $\bb{V} \subseteq \bb{S5}_2$, let $\bb{V}_\text{SI}$ be the class of s.i.~$\bb{V}$-algebras, and define $T(\bb{V})$ to be the variety generated by the translations of the s.i.~members:
\[
T(\bb{V}) = HSP(\setbuilder{T(\frk{A})}{\frk{A} \in \bb{V}_\text{SI}}).    
\]
It follows from \cref{lem:translation-descriptive-frame} that $T(\bb{V})$ is generated by a class of $\msfl[2]$-algebras, so $T(\bb{V}) \subseteq \msfv[2]$.

\begin{theorem}
    \label{thm:translation-local-finiteness}
    $T$ preserves and reflects local finiteness; that is, $\bb{V} \subseteq \bb{S5}_2$ is locally finite iff $T(\bb{V}) \subseteq \msfv[2]$ is locally finite.
\end{theorem}

\begin{proof}
    Suppose $T(\bb{V})$ is locally finite, and let $\frk{A} = \frk{F}^* \in \bb{V}_\text{SI}$ be generated by $g_1, \dots, g_n$.
    Then each $g_i$ is 
    a subset of $D_2 \subseteq T(\frk{F})$, and by \cref{lem:second-layer} the operations of $\frk{A}_2 \cong \frk{A}$ are definable using the clopen set $D_2$.
    By assumption, there is $f:\omega\to\omega$ such
    that the size of the subalgebra $\la g_1, \dots, g_n, D_2 \ra \subseteq T(\frk{A})$ is bounded by $f(n+1)$. This subalgebra contains
    every element of
    $\frk{A}$. Thus, $\abs{\frk{A}} \leq f(n+1)$.
    By \cite[Thm.~3.7(4)]{Bez01}, 
    we conclude that $\bb{V}$ is locally finite.

    Suppose $\bb{V}$ is locally finite. To prove that $T(\bb{V})$ is locally finite, it is sufficient to 
    show that the generating class $\setbuilder{T(\frk{A})}{\frk{A} \in \bb{V}_\text{SI}}$
    is uniformly locally finite (see \cite[285]{Mal73}). 
    Suppose $\frk{B}$ is an $n$-generated subalgebra of $T(\frk{A})$, corresponding to a correct partition $L$ of $T(\frk{A}_*)$.
    By \cref{lem:loc-fin-right-direction}, there is a $2n$-generated subalgebra $\frk{B}' \subseteq \frk{A}$ such that the corresponding correct partition $K$ of $\frk{A}_*$ satisfies $\what{K} \subseteq L$. Therefore, $\frk{B} \subseteq T(\frk{B}')$. By assumption, we have $f:\omega\to\omega$ such
    that $\abs{\frk{B}'} \leq f(2n)$.
    Thus, it follows from Construction 6.8 that the size of $T(\frk{B'})$ and hence of $\frk{B}$ is bounded by $f(2n)^2$.
\end{proof}

\section{Conclusion and future work}\label{sec: conclusion}

We have shown that, though the Segerberg--Maksimova theorem generalizes to certain subvarieties of $\bb{MS4}$, it fails drastically for $\bb{MS4}$ in general, and already in $\msfv$.
As shown in \cref{translation}, characterizing local finiteness even in $\msfv[2]$ is at least as hard as in $\bb{S5}_2$.
It is then natural to ask if the problem of local finiteness in $\bb{MS4}$ is \textit{strictly} harder; i.e., would a characterization of local finiteness in $\bb{S5}_2$ yield one for $\bb{MS4}$?
We finish with an example that suggests a negative answer.

The frame $T(\frk{F})$ constructed in \cref{translation} has the property that the relativization $\frk{A}_1$ to $D_1$ is a locally finite $\msf{S5}^2$-algebra, while the relativization $\frk{A}_2$ to $D_2$ is the original $\msf{S5}_2$-algebra $\frk{F}^*$.
In light of \cref{thm:translation-local-finiteness}, the algebra $T(\frk{F})^*$ is locally finite iff $\frk{A}_2$ is locally finite, which happens iff the relativization to each of the layers of $T(\frk{F})$ is locally finite.
It is then natural to conjecture that the local finiteness of an $\msf{MS4}$-algebra (of finite depth) is completely determined by the local finiteness of the $\msf{S5}_2$-algebras of the relativizations to each layer.
We show this is not the case already in $\msfv[3]$.
Indeed, the following figure gives an example of an $\msfl[3]$-frame $\frk{F}$ whose relativization to each of the three layers is a locally finite $\msf{S5}_2$-algebra, and yet $\frk{F}^*$ fails to be locally finite:

\begin{figure}[h]
\begin{tikzpicture}[
    scale=1.25,
    dot/.style={circle,fill=black,minimum size=6pt,inner sep=0}
]

\foreach \x in {0,...,4,6} {
    \draw[blue,thick,fill=blue!20] (\x - 0.4, -1.25) -- (\x - 0.4, 1.85) -- (\x + 0.15, 2.85) -- (\x + 0.15, -0.25) --cycle;
}

\foreach \y in {2.05, 1, -0.05} {
    \draw[dotted] (-0.6, \y + 0.2) -- (6.6, \y + 0.2) -- (6.2, \y - 0.75) -- (-1, \y - 0.75) -- cycle;
}

\node (C2) [dot] at (-0.2, 0.4) {};
\node (D2) [dot] at (0.8, 0.4) {};
\node (E2) [dot] at (1, 1) {};
\node (F3) [dot] at (2, 0) {};
\node (G2) [dot] at (1.8, 0.4) {};
\node (H2) [dot] at (2.8, 0.4) {};
\node (I2) [dot] at (3, 1) {};
\node (J3) [dot] at (4, 0) {};
\node (K2) [dot] at (3.8, 0.4) {};
\node at (5, 0.75) {$\dots$};
\node at (5, -0.25) {$\dots$};
\node (LIM2) [dot] at (6, 1) {};
\node (LIM3) [dot] at (6, 0) {};

\node (A1) [dot] at (0, 2) {};
\node (B1) [dot] at (1, 2) {};
\node (C1) [dot] at (2, 2) {};
\node (D1) [dot] at (3, 2) {};
\node (E1) [dot] at (4, 2) {};
\node at (5, 2) {$\dots$};
\node (LIM1) [dot] at (6, 2) {};

\draw (C2) -- (D2);
\draw[-latex] (F3) -- (E2);
\draw[-latex] (J3) -- (I2);
\draw (G2) -- (H2);

\draw (A1) -- (E1);
\draw (E1) -- (4.5, 2) {};
\draw (5.5, 2) -- (LIM1) {};

\draw[-latex] (E2) -- (B1);

\draw[-latex] (C2) to[out=90, in=250] (A1);
\draw[-latex] (D2) to[out=90, in=250] (B1);
\draw[-latex] (G2) to[out=90, in=250] (C1);
\draw[-latex] (H2) to[out=90, in=250] (D1);
\draw[-latex] (K2) to[out=90, in=250] (E1);
\draw[-latex] (LIM2) to (LIM1);
\draw[-latex] (LIM3) to (LIM2);

\draw (C2) circle[radius=0.2];
\node at (-0.6, 0.4) {$g$};

\end{tikzpicture}
\caption{An $\msfl[3]$-frame $\frk{F}$ in which the relativization to each layer is a locally finite $\msf{S5}_2$-algebra, but $\frk{F}^*$ is not locally finite}
\label{fig:3d-snake}
\end{figure}

\begin{example}
    Let $\frk{F}$ be the frame in \cref{fig:3d-snake}.
    Topologically, $\frk{F}$ is the disjoint union of three layers, where the topology on each layer makes the right-most point the limit point of a one-point compactification of the finite-index points with the discrete topology.
    An argument similar to the proof of \cref{lem:translation-descriptive-frame} shows that $R$ and $E$ are continuous relations,
    and the ``top rail'' ensures that $\frk{F}$ is an $\msf{MS4}$-frame. Clearly the $R$-depth of $\frk{F}$ is $3$.
    
    Let $g$ be the singleton set depicted in the figure and $d$ the top layer.
    Then we can generate an infinite family of sets from the two clopen sets $g$ and $d$ as follows:
    \[
        s_0 = g \qquad s_{n+1} = \exists \lozenge s_n - d.
    \]
    Thus, $\frk{F}^*$ is not locally finite.
    
    However, the relativization $\frk{A}_i$ to each layer is easily seen to be a locally finite $\msf{S5}_2$-algebra.
    As remarked in the discussion preceding \cref{lem:second-layer}, $\frk{A}_1$ belongs to a proper subvariety of $\bb{S5}^2$, hence is locally finite.
    The second and third layers, equipped with the restrictions of the relations $R$ and $E$, form an $\msf{S5}_2$-frame that is a compactification of an infinite disjoint union of a single finite frame (a 3-element frame in the case of the second layer, and a singleton frame in the case of the third).
    Thus, for $i \in \set{2,3}$, $\frk{A}_i$ is a subalgebra of an infinite power of a single finite algebra, so $\frk{A}_i$ belongs to a finitely generated variety, 
    and hence $\frk{A}_i$ is locally finite.
\end{example}

There are two natural directions of research suggested by the results in this paper.
First, as was demonstrated in \cite{BBI23}, the monadic version of Casari's formula plays an important role in obtaining a faithful provability interpretation of the one-variable fragment of intuitionistic predicate logic. Let $\msf{M^{+}S4}$ be the extension of $\msf{MS4}$ obtained by postulating the G\"odel translation of this formula. The restrictions on the interaction of the relations $R$ and $E$ in $\msf{M^{+}S4}$-frames make it plausible to expect a more manageable characterization of locally finite varieties of $\msf{M^{+}S4}$-algebras akin the Segerberg--Maksimova theorem.

Another direction of research is suggested by the nature of the non-locally finite frames given in the paper.
The $\msf{S5}_2$-frame in \cref{fig:translation} and the frame in \cref{fig:3d-snake} both exhibit a one-generated infinite subalgebra in a similar manner, by alternating through the (non-commuting) relations.
It is worth investigating whether this phenomenon is the only reason for non-local finiteness in $\bb{S5}_2$, and even in $\bb{MS4}$. 

\section*{Acknowledgment}

We thank the referee for careful reading and useful feedback. 

\printbibliography

\end{document}